\documentclass{amsart}
\usepackage{amsmath,amsfonts,amssymb,amscd,amsthm}
\newtheorem{theorem}[equation]{Theorem}
\newtheorem{lemma}[equation]{Lemma}
\newtheorem{proposition}[equation]{Proposition}
\newtheorem{definition}[equation]{Definition}
\newtheorem{corollary}[equation]{Corollary}
\theoremstyle{remark}
\newtheorem{remark}[equation]{Remark}
\numberwithin{equation}{section}

\DeclareMathOperator{\Id}{Id}
\DeclareMathOperator{\Ind}{Ind}
\DeclareMathOperator{\Res}{Res}
\DeclareMathOperator{\Int}{Int}
\DeclareMathOperator{\Irr}{Irr}
\DeclareMathOperator{\Trace}{Trace}
\DeclareMathOperator{\codim}{codim}
\DeclareMathOperator{\class}{class}

\newcommand{\bbF}{{\mathbb F}}
\newcommand{\bbR}{{\mathbb R}}
\newcommand{\bbC}{{\mathbb C}}
\newcommand{\bbQ}{{\mathbb Q}}
\newcommand{\bbZ}{{\mathbb Z}}
\newcommand{\bG}{{\bf G}}
\newcommand{\bH}{{\bf H}}
\newcommand{\bL}{{\bf L}}
\newcommand{\bM}{{\bf M}}
\newcommand{\bT}{{\bf T}}
\newcommand{\bU}{{\bf U}}
\newcommand{\ci}{{\iota}}
\newcommand{\cj}{{\gamma}}
\newcommand{\ck}{{\kappa}}
\newcommand{\cregG}{{\rho^\bG}}
\newcommand{\cregM}{{\rho^\bM}}
\newcommand{\creg}{{\rho}}
\newcommand{\cs}{{\sigma}}
\newcommand{\CC}{{\mathcal C}}
\newcommand{\cD}{{\mathcal D}}
\newcommand{\CF}{{\mathcal F}}
\newcommand{\CI}{{\mathcal I}}
\newcommand{\CP}{{\mathcal P}}
\newcommand{\CX}{{\mathcal X}}
\newcommand{\CY}{{\mathcal Y}}
\newcommand{\CIM}{{\CI_\bM}}
\newcommand{\CJG}{{\CC_\CI(\GF)}}
\newcommand{\CuG}{{\Cuni(\GF)}}
\newcommand{\Cuni}{{\CC_{\text{uni}}}}
\newcommand{\CWG}{\CC(\WGLF)}
\newcommand{\Fq}{{\bbF_q}}
\newcommand{\GF}{{\bG^F}}
\newcommand{\MF}{{\bM^F}}
\newcommand{\Mo}{{\bM_0}}
\newcommand{\Qlbar}{{\overline{\mathbb Q}_\ell}}
\newcommand{\semi}{{\rtimes}}
\newcommand{\sgn}{{\varepsilon}}
\newcommand{\GL}{{\text{GL}}}
\newcommand{\SL}{{\text{SL}}}
\newcommand{\sR}{{\lexp *R}}
\newcommand{\WGLF}{{\WGL.F}}
\newcommand{\WGL}{{W_\bG(\bL)}}
\newcommand{\WMLF}{{\WML.wF}}
\newcommand{\WML}{{W_\Mo(\bL)}}
\newcommand{\ZL}{{{\mathcal Z}_\bL}}
\newcommand{\gothS}{{\mathfrak S}}
\newcommand{\inv}{^{-1}}
\newcommand{\cf}{{\it cf.}}
\newcommand{\eg}{{\it e.g.}}
\newcommand{\ie}{{\it i.e.}}
\newcommand{\lexp}[2]{\kern\scriptspace\vphantom{#2}^{#1}\kern-\scriptspace#2}
\newcommand{\scal}[3]{{\langle\,#1,#2\,\rangle_{#3}}}
\newcommand{\te}{{\tilde\sgn}}
\newcommand{\tG}{{\tilde\Gamma}}
\newcommand{\tX}{{\tilde\CX}}
\newcommand{\tphi}{{\tilde\varphi}}
\newcommand{\tQ}{{\tilde Q}}
\newcommand{\genby}[1]{\mathopen<#1\mathclose>}
\begin{document}
\author{F.~Digne \and G.~Lehrer \and J.~Michel}
\title[The space of unipotently supported class functions]{The space of unipotently supported class functions
on a finite reductive group}
\date{May 31, 2002}
\dedicatory{to Robert Steinberg}
\maketitle
\section{Introduction}
Let $\bG$  be a  connected reductive algebraic  group over  an algebraic
closure  $\bbF$  of the  finite  field  $\Fq$  of $q$  elements;  assume
$\bG$  has an  $\Fq$-structure  with  associated Frobenius  endomorphism
$F$  and let  $\ell$  be a  prime distinct  from  the characteristic  of
$\Fq$.  In \cite[\S  7.1]{DLM1} and  \cite{DLM2} we  outlined a  program
for  the determination  of  the irreducible  $\Qlbar$-characters of  the
finite  group  $\GF$, which  showed  that  the  problem may  be  largely
reduced  (by induction)  to  an explicit  determination  of the  Lusztig
restrictions  $\sR^\bG_\bM(\chi)$  of  all  the  irreducible  characters
$\chi$ of $\GF$,  for all rational Levi subgroups $\bM$  of $\bG$. Here,
and throughout  this paper, the  word ``rational'' means  ``stable under
the  action of  $F$''.  As shown  in \cite{DLM2},  this  problem may  be
addressed  through the  determination of  the Lusztig  restrictions $\sR
^\bG_\bM(\Gamma_u)$, where  $\Gamma_u$ is the  generalized Gelfand-Graev
character  corresponding to  the $\GF$-conjugacy  class of  the rational
unipotent element $u\in \GF$.

Now the characters  $\Gamma_u$ are examples of class  functions on $\GF$
which vanish  outside the  unipotent set. Such  functions form  a vector
space  over $\Qlbar$,  which we  denote by  $\CuG$; it  is the  space of
unipotently  supported class  functions  on $\bG$.  The $\Gamma_u$  form
a  basis  of  this  space,  and  our  strategy  in  this  work  will  be
to  determine  the  map $\sR_\bM^\bG:\CuG\to\Cuni(\MF)$  explicitly.  We
shall  use Lusztig's orthogonal decomposition of  the space  $\CuG$ into  summands
corresponding to  ``rational blocks'' (see below) and  determine $\sR_\bM^\bG$ on
each block generically,  \ie\ in terms of Weyl group  data which is 
associated with the block. In  particular, we obtain a  simple expression
for  the Lusztig  restriction of  generalized Green  functions. We  then
express the generalized Gelfand-Graev characters  in terms of this basis
to  describe  their  Lusztig  restriction. In  \cite{DLM2}  we  computed
$\sR_\bM^\bG$ of the generalized Gelfand-Graev character which corresponds
to a regular unipotent class.  In this work, we apply the  general method to carry
out the corresponding computation explicitly in the subregular case.

Our  general  result  on   $\sR_\bM^\bG$  of  generalized  Gelfand-Graev
characters  (\ref{*RMG(Gamma)})  essentially  reduces  this  computation
to  the  two problems  of  finding  the Poincar\'e  polynomials  $\tilde
P_{\ci,\ck}$ of certain intersection  complexes on closures of unipotent
classes,  and to  the  computation of  induction-restriction tables  for
twisted  characters of  Weyl groups.  In  \S 8  we also  prove a  result
(\ref{Sln}  below)  which reduces  these  computations  in the  case  of
$\SL_n$ to  the case of  $\GL_{n'}$, for various $n'$.  These investigations
are part of our strategy of reducing the computation of character values
to the case of ``high'' unipotent classes in the usual partial order.

The first five sections of this paper consist largely of a recasting
of the of work of
Lusztig, which may be found in \cite{US},\cite{ICC},\cite{CS}, in a form
which  permits  practical computation. They also contain several orthogonality
relations for Green functions and their generalizations, which are proved by relating
the inner product in $\CuG$ to the inner product of twisted class functions
on a Weyl group. In section 6, we prove orthogonality relations for
the generalized Gelfand-Graev characters in the same way, in addition to
determining their Lusztig restriction. By  and  large we  maintain  the
notation of  \cite{DLM2}. We  shall rely on  the context  to distinguish
between the Frobenius  endomorphism $F$ of an $\Fq$-group  $\bG$ and the
automorphisms, also denoted $F$, which  are induced by $F$ on reflection
groups  (such  as the  Weyl  group)  which  are associated  with  $\bG$.
Throughout this work we shall freely  use the character theory of cosets
of a finite group, for which the reader is referred to \cite[(0.4)]{DM2}
or  \cite{Ka1}.  Characters  of  cosets  are  also  sometimes  known  as
``twisted class functions''.

\section{Preliminaries}
Let $\ci=(C,\zeta)$ be  a pair consisting of a unipotent  class of $\bG$
and an  irreducible $\bG$-equivariant  $\Qlbar$-local system  $\zeta$ on
it; then  $C$ will be  called the {\it  support} of $\ci$  and sometimes
denoted $C_\ci$. If we fix  a non-trivial additive character $\chi_0$ of
the  prime  field $\bbF_p$  of  $\Fq$,  as  in \cite[1.6]{DLM2}  we  may
define a generalized Gelfand-Graev function $\Gamma_\ci$ associated with
$\ci$; one  of our  objectives here is  to express  Lusztig restrictions
of  generalized   Gelfand-Graev  characters  in  terms   of  generalized
Gelfand-Graev characters.

As in  \cite{DLM2}, if  the pair  $\ci$ is  $F$-stable, we  shall follow
Lusztig  \cite[24.1--24.2]{CS}  in  making   a  specific  choice  of  an
isomorphism  $\sigma:F^*\zeta\xrightarrow\sim\zeta$,  and we  denote  by
$\CY_\ci$  the  characteristic  function of  $\zeta$  which  corresponds
to  $\sigma$,  and  by  $\CX_\ci$ the  characteristic  function  of  the
intersection  cohomology   complex  of  $\zeta$  (for   $u\in  C^F$,  we
have  $\CX_\ci(u)=\CY_\ci(u)$). The  set  $\CP$ of  all  pairs $\ci$  is
partitioned  into ``blocks''  $\CI$,  each of  which  has an  associated
cuspidal  datum  $(\bL,\ci_0=(C_0,\zeta_0))$  where   $\bL$  is  a  Levi
subgroup of  some parabolic  subgroup of  $\bG$, which  is unique  up to
$\bG$-conjugacy. If the  block concerned is rational,  then as explained
in  \cite[1.4]{DLM2},  both $\bL$  and  the  parabolic subgroup  may  be
assumed to  be rational. The pairs  in the block $\CI$  are in bijection
with  the irreducible  characters  of  the group  $\WGL=N_\bG(\bL)/\bL$,
which is a Coxeter group. If $\CI$ is a rational block and $\varphi_\ci$
is  the character  associated  in this  way to  $\ci=(C,\zeta)\in\CI^F$,
then  an  extension  $\tphi_\ci$ of  $\varphi_\ci$  to  $\WGL\semi\genby
F$ determines  an isomorphism  $\sigma:F^*\zeta\xrightarrow\sim\zeta$ as
above.  In this  work,  we shall  always choose  $\tphi_\ci$  to be  the
``preferred extension'' described in \cite[17.2]{CS} (as Lusztig does in
\cite[24.2]{CS}).

The  functions  $\CY_\ci$ form  a  basis  of  the space  of  unipotently
supported  class  functions  on  $\GF$   as  $\ci$  runs  over  the  set
$\CP^F$  of   all  rational  pairs.   For  a  given  block   $\CI$,  the
functions  $\CX_\ci$   form  another   basis  of  the   space  spanned
by  $\{\CY_\ci\}_{\ci\in\CI^F}$,  and   if  we  write  $\CX_\ci=\sum_\ck
P_{\ck,\ci}\CY_\ck$    then    the   $P_{\ck,\ci}$    are    polynomials
in   $q$   with   integer    coefficients   We   have   $P_{\ck,\ci}=0$
unless   $C_\ck\subset\overline  C_\ci$   and   if  $C_\ck=C_\ci$   then
$P_{\ck,\ci}=\delta_{\ck,\ci}$ (see \eg, \cite[6.5]{US}). We will assume
from now on that the pairs $\ci$ have been totally ordered in such a way
that  $C_\ck\subset\overline  C_\ci  \Rightarrow  \ck\le\ci$.  Then  the
matrix $(P_{\ck,\ci})$ is upper unitriangular.

Set $\tX_\ci=q^{c_\ci}\CX_\ci$ and $\tilde\CY_\ci=q^{c_\ci}\CY_\ci$
where $c_\ci=\frac12(\codim C_\ci-\dim Z_\bL)$.
Then we have $\tX_\ci=\sum_\ck \tilde P_{\ck,\ci}\tilde\CY_\ck$,
where $\tilde P_{\ck,\ci}=q^{c_\ci-c_\ck} P_{\ck,\ci}$.

\begin{remark}\label{innerprod}
We shall speak below of ``complex conjugation'' in the field $\Qlbar$,
denoted by   $a\mapsto\overline a$. This is justified by noting that 
$\Qlbar$ is abstractly isomorphic to $\bbC$. In practice, we shall
apply this notion almost exclusively to the subfield of $\Qlbar$ which is 
generated by all roots of unity, on which conjugation is uniquely defined since it
fixes $\bbQ$ and inverts roots of unity. 
 We  therefore speak of ``real''  values (meaning
fixed by conjugation) and ``complex conjugates'' in this context.
The  space  $\CuG$  is then an  inner  product  space  with  Hermitian
form  defined by  $$\scal{f}{g}\GF=|\bG|\inv\sum_{x\in \GF}f(x)\overline
{g(x)}.$$\end{remark}

\begin{remark}\label{twist}
The cuspidal  datum $(\bL,\ci_0)$ defines  a unique block $\CIM$  of any
Levi subgroup $\bM$ of $\bG$  which contains a $\bG$-conjugate of $\bL$.
Assume $\bM$  rational, and let  $\bL'=\Int g(\bL) (:=g\bL g\inv)$  be a
conjugate  of  $\bL$ which  is  rational  and  contained in  $\bM$;  let
$\Mo\supset\bL$  be the  conjugate  $\Int g\inv(\bM)$  of $\bM$.  Define
$w\in \WGL$  by $\dot  w=g\inv F(g)\in  N_\bG(\bL)$. Then  $(\bL',F)$ is
conjugate to $(\bL,\dot wF)$ and  $\Mo$ is $\dot wF$-stable; moreover we
may identify (via $\Int g\inv$) $(\bM,F)$ with $(\Mo,\dot wF)$ and hence
$(W_\bM(\bL'),F)$ with  $(\WML,wF)$, $\Cuni(\MF)$  with $\Cuni(\Mo^{\dot
wF})$  and $(\CIM,F)$  with $(\CI_\Mo,wF)$.  A particular  case of  this
occurs when $\Mo=\bL$, when we refer  to the twisted version of $\bL$ as
$\bL_w$  (for  $w\in \WGL$).  The  cuspidal  pair  $\ci_0$ of  $\bL$  is
taken  by $\Int(g)$  to a  cuspidal pair  of $\bL_w$.  The corresponding
characteristic function on  $\bL_w^F$ is likewise taken by  $g\inv$ to a
function on $\bL^{\dot wF}$, which we denote by $\CX_{\ci_0,w}$.
\end{remark}

We recall that Lusztig induction  $R^\bG_\bM$ has an easy description in
terms of the  functions $\CX_\ci$, which applies  with some restrictions
on $p$ and $q$. The results of this paper will depend on this, and hence
we shall  assume, sometimes without  explicit mention, for the  whole of
our work that (\cf\ \cite[3.1]{DLM2})  the characteristic $p$ is good for
$\bG$ and that $q>q_0(\bG)$, a constant which depends only on the Dynkin
diagram of $\bG$.

\begin{proposition}\label{a} Assume $p$ good and  $q$ sufficiently
large, and that $\bM$ contains a rational conjugate $\bL_w$ of $\bL$
as in \ref{twist}. Assume (as we may, by the above discussion)
that $\bL_w$ is a split Levi subgroup of $\bM$. 
Then for $\ci\in\CIM^F$, we have:
\begin{enumerate}
\item$ \displaystyle R_\bM^\bG(\tX_\ci)=\sum_{\ck\in\CI^F}
\scal{\tphi_\ci}{\Res^\WGLF_\WMLF\tphi_\ck}\WMLF\tX_\ck,$
where $R_\bM^\bG$ is the Lusztig induction functor.
\item
$\scal{\tphi_\ci}{\Res^\WGLF_\WMLF\tphi_\ck}\WMLF=0$
unless $\overline{C_\ci}\subset\overline{C_\ck}\subset
\overline{\Ind_\bM^\bG C_\ci}$.
\end{enumerate}
\end{proposition}
\begin{proof} Assertion (i) is in  \cite[3.3]{DLM2}. Let us prove (ii).
For  the rightmost  inclusion recall  that, from  the definition  of the
induction of  perverse sheaves,  only pairs  $\ck$ with  support smaller
than that of  the class induced from  the support of $\ci$  can have non
zero  coefficient  in  $R_\bM^\bG(\tX_\ci)$. To  prove  the  other
inclusion, first notice that if
$\scal{\tphi_\ci}{\Res^\WGLF_\WMLF\tphi_\ck}\WMLF$
is non-zero then so is
$\scal{\varphi_\ci}{\Res^\WGL_\WML\varphi_\ck}\WML$.
But it follows from formula (II)  in \cite[1.2]{S} that the latter inner
product  is zero  unless there  exists  a representative  of $C_\ck$  in
$C_\ci.\bU$  where  $\bU$  is  the  unipotent  radical  of  a  parabolic
subgroup  admitting $\bM$  as a  Levi  component. This  in turn  implies
$C_\ci\subset\overline{C_\ck}$ by \cite[5.8]{DLM1}.
\end{proof}
\begin{remark}\label{split}
We shall often have a situation where $\bM$ is a rational Levi subgroup
of $\bG$ which contains a rational conjugate $\bL_w$ of $\bL$, as in 
\ref{twist}. In this situation we shall consistently assume $w\in \WGL$
to have been chosen so that $\bL_w$ is split in $\bM$, i.e., is contained 
in a rational parabolic subgroup of $\bM$. In this case 
$w\in \WGL$ is determined up to $F$-conjugacy in $\WGL$ and the function
$R_{\bL_w}^\bG(\CX_{\ci_0,w})$ is well defined (see \cite[3.2 and 
3.3(1)]{DLM2}). This is implicit in the statement and proof of \ref{a}.
\end{remark}

\section{Generalized Green functions and Lusztig restriction}
In this section we shall  interpret Lusztig induction and restriction in
terms ordinary induction  and restriction of twisted  class functions on
cosets of  parabolic subgroups of Coxeter  groups. This will be  done by
defining  a  linear isomorphism  between  the  spaces of  twisted  class
functions on $\WGL$  and a certain subspace of the  space of unipotently
supported  functions. Under  this map,  the (normalized)  characteristic
functions of the $F$-classes of $\WGL$ correspond to functions we define
as ``generalized Green functions''. These  are analogues of the ordinary
Green functions  (the latter  corresponding to the  ``principal block'')
which constitute  a basis  of the space  of unipotently  supported class
functions.  In order  to  compute their  Lusztig  restriction, we  shall
relate the generalized Gelfand-Graev characters to these.

For  the  whole of  this  section,  we  fix  a rational  cuspidal  datum
$(\bL,\ci_0)$,  where  we  may  assume  that $\bL$  is  split,  \ie\  is
contained in  a rational parabolic  subgroup of $\bG$.  Let $\CC(\WGLF)$
be  the  space  of  $\WGL$-invariant functions  (\ie\  class  functions)
on  $\WGLF$  and  recall  that   $\CuG$  is  the  space  of  unipotently
supported  class functions  on $\GF$.  For each  $w\in \WGL$,  we fix  a
$w$-twisted rational conjugate  $\bL_w$ of $\bL$ as  in \ref{twist}, \ref{split}, and
$\tX_{\ci_0,w}\in\Cuni(\bL_w^F)$ is the  class function on $\bL^F$
(see \ref{twist} and \ref{split}) associated with $\ci_0$.

\begin{definition}\label{defqg} Let $\CJG$
be the subspace of  $\CuG$  spanned  by  the  functions
$\{\CY_\ci\mid \ci\in\CI^F\}$.
\begin{enumerate}
\item   Define  the linear  isomorphism  $Q^\bG$  from   $\CC(\WGLF)$  to
$\CJG$  by $Q^\bG(\tphi_\ci)=\tX_\ci$.
\item  For  $w\in\WGL$  define  $\gamma_{wF}\in\CC(\WGLF)$  by
$$\gamma_{wF}(vF)=\begin{cases}
0,&\text{if $vF$ is  not $\WGL$-conjugate to $wF$}\\
|C_\WGL(wF)|,&\text{otherwise.}
\end{cases}$$
\item The {\bf  generalized Green
function} $Q^\bG_{wF}$ is defined by $Q^\bG_{wF}=Q^\bG(\gamma_{wF})$.
\end{enumerate}
\end{definition}
Note  that   since  the  (distinct)   $\gamma_{wF}$  form  a   basis  of
$\CC(\WGLF)$, the generalized Green  functions $Q^\bG_{wF}$ form a basis
of $\CJG$.

We shall omit the superscript in  $Q^\bG$ and $Q_{wF}^\bG$ when there is
no ambiguity.

\begin{proposition}\label{defQwF}
We have $Q_{wF}=R_{\bL_w}^\bG \tX_{\ci_0,w}$.
\end{proposition}
\begin{proof}
Since the  $\tphi_\ci$ form an orthonormal basis of $\CC(\WGLF)$,
and $\scal\theta{\gamma_{wF}}\WGLF=\theta(wF)$ for any
$\theta\in \CC(\WGLF)$, we have
$$\gamma_{wF}=\sum_{\ci\in\CI^F}\scal{\tphi_\ci}{\gamma_{wF}}\WGLF\tphi_\ci=
\sum_{\ci\in\CI^F}\tphi_\ci(wF)\tphi_\ci,$$ whence by linearity
\begin{equation}\label{QwF}
Q_{wF}=\sum_{\ci\in\CI^F}\tphi_\ci(wF)\tX_\ci.
\end{equation}

But by  \cite[3.1]{DLM2} we have
\begin{equation}\label{CXi}
\tX_\ci=|\WGL|\inv\sum_{v\in\WGL}\tphi_\ci(vF)R_{\bL_v}^\bG(\tX_{\ci_0,v}).
\end{equation}
Now  in \ref{CXi},  the  summand corresponding  to  $w\in \WGL$  depends
only  on  the $\WGL$-class  of  $wF$.  To  see  this, observe  that  the
function  $\tX_{\ci_0,v}$   is  invariant  under   conjugation  by
$N_\bG(\bL_v)^F$,  so that  $R_{\bL_v}^\bG \tX_{\ci_0,v}$  depends
only  on the  $\GF$-class  of  $\bL_v$, which  is  parametrized by  the
$W$-class of the coset $W_\bL.vF$, or by the $\WGL$-class of the element
$vF\in \WGLF$.

Since the $\tphi_\ci$ take real values, the second orthogonality
relation for them reads
$$\sum_\ci\tphi_\ci(wF) \tphi_\ci(vF)=
\begin{cases}0,&\text{ if $vF$ is  not $\WGL$-conjugate to $wF$}\\
|C_\WGL(wF)|,&\text{otherwise.}
\end{cases}$$
Substituting  \ref{CXi}  into \ref{QwF}  and  using  this relation,  the
result follows.
\end{proof}
It   follows  from   this   proposition  that   our  generalized   Green
functions   are   the  same   as   those   in  \cite[8.3.1]{CS},   since
$q^{c_{\ci_0}}\CX_{\ci_0}$ is the restriction  to the unipotent elements
of   the  characteristic   function  of   the  perverse   sheaf  denoted
by  $\mathop{\rm   IC}(\overline\Sigma,{\mathcal  E})[\dim(\Sigma)]$  in
\cite[8.2]{CS}  and  for  cuspidal local  systems,  Lusztig's  induction
coincides with the induction of perverse sheaves by \cite{GF}.

Both $\CWG$ and  $\CJG$ have natural structures  as non-degenerate inner
product  spaces. Although  $Q^\bG$ is  not  an isometry,  its effect  on
scalar products can be computed.

\begin{definition}\label{defCZGL}
Define the function $\ZL\in \CWG$ by $\ZL(wF)=|Z^{0wF}_\bL|=|Z_{\bL_w}^{0F}|$.
\end{definition}

\begin{proposition}\label{scalar product QG} We have, for any
two functions $\theta,\phi\in \CWG$,
$$\scal{Q^\bG(\theta)}{Q^\bG(\phi)}\GF=\scal{\ZL\inv\theta }\phi\WGLF.$$
\end{proposition}
\begin{proof}
First note that \cite[24.3.6]{CS}, suitably interpreted
to take into account the distinction between our $\tX_{\ci_0}$
and Lusztig's $\CX_{\ci_0}$, shows that
\begin{equation}\label{scalX}
\scal{\tX_\ci}{\tX_\ck}\GF=\scal{\ZL\inv\tphi_\ci}{\tphi_\ck}\WGLF.
\end{equation}
Now in order to prove the proposition,  it suffices to do so as $\theta$
and $\phi$  run over a  basis of $\CWG$  In particular, it  suffices to
take $\theta=\tphi_\ci$  and $\phi=\tphi_\ck$.  But then
the statement is precisely the equation \ref{scalX}, whence the result.
\end{proof}
It follows  easily from the definition  \ref{defqg}(iii) and
\ref{scalar product QG} that  the generalized  Green functions  form an
orthogonal basis of $\CJG$. More precisely, we have
\begin{corollary}\label{scalar product Q}
$$\scal{Q_{wF}}{Q_{w'F}}\GF=
 \begin{cases}0&\text{if $wF$ and $w'F$ are not conjugate in $\WGL$,}\\
 \dfrac{|C_\WGL(wF)|}{|Z^{0wF}_\bL|}&
 \text{otherwise.}
 \end{cases}$$
\end{corollary}

The formula  \ref{scalar product  Q} superficially seems  different from
\cite[9.11]{CS}.  However  the  two formulae  are  actually  equivalent,
although there is a power of $q$ in {\it loc.cit.} which is absent here.
This is explained by the facts  that in [{\it loc.cit.}, 9.11] the inner
product used differs from ours, in that it does not involve conjugation,
and  that the  formula  given there  is  for the  inner  product of  two
Green  functions corresponding  to  contragredient  local systems,  with
contragredient  Frobenius isomorphisms.  In Lusztig's  notation, if  the
characteristic function  of the  sheaf $\CF$ with  Frobenius isomorphism
$\varphi_1$  is  $f$, then  the  characteristic  function of  $\CF^\vee$
with Frobenius isomorphism $\varphi_1^\vee$ is $q^{-2c_{\ci_0}}\overline
f$  (see the  computation  in  the proof  of  \cite[9.8]{CS}); this,  in
conjunction  with  the  fact  that  $R^\bG_\bM$  commutes  with  complex
conjugation, shows the formulae are equivalent.

\begin{remark}\label{epsilon}  The  preferred  extension  $\te$  of  the
alternating  character $\sgn$  of $\WGL$  will play  a prominent  r\^ole
in  our  work.  A  fact  which   we  shall  use  repeatedly,  and  which
results  from  the  description  in  \cite[17.2]{CS}  of  the  preferred
extension,  is that  $\te$  is  trivial on  Frobenius,  i.e., for  $w\in
\WGL$, $\te(w.F)=\sgn(w)$. Note also  that since the preferred extension
is  real,  if  $\tphi_\ci$  is  the  preferred  extension  corresponding
to  $\ci\in\CI^F$, then  there  is  a sign  $\sgn_\ci=\pm  1$ such  that
$\tphi_\ci\otimes   \te=\sgn_\ci\tphi_{\hat\ci}$,  where   $\hat\ci$  is
defined by $\varphi_\ci\otimes\sgn=\varphi_{\hat\ci}$.
\end{remark}
\smallskip

Let $\bH$ be any linear algebraic group with a Frobenius morphism
$F:\bH\to\bH$ which corresponds to an $\Fq$-structure on $\bG$. Let $\bT$
be a maximally split maximal torus of $\bH$ and write $R_u(\bH)$ for the
unipotent radical of $\bH$. Then the Weyl group $W=W_\bH(\bT)$ acts
as a reflection group on 
$Y(\bT)\otimes\bbR$, and $F$ has an induced action as $q\phi$ on this space,
where $\phi$ is a linear transformation of finite order 
(\cf\ \cite[p.40]{DM1}). Write
$\{f_1,f_2\dots,f_\ell\}$ for a set of basic invariants of $W$ and
let $d_i=\deg(f_i)$. It is known (\cf\ \cite[6.1]{Sp}) that the $f_i$
may be chosen to be eigenfunctions for $\phi$, i.e. $\phi f_i=\delta_i f_i$
for each $i$, where $\delta_i\in \bbC$.

\begin{lemma}\label{orderinv} With notation as in the previous paragraph,
we have 
\begin{enumerate}
\item The order of $\bH$ is given by 
$$
|\bH^F|=q^{\dim R_u(\bH)+\sum_i(d_i-1)}\prod_i(q^{d_i}-\delta_i).
$$
\item
If $F$ is varied by keeping $\phi$ fixed and allowing $q$ to vary, the
order function in (i) is a polynomial in $q$ and 
$$
|\bH^F|(q^{-1})=q^{-\dim \bH}\sgn_\bH |\bH^F|_{q'},
$$
where, for any linear algebraic group $\bH$ we write
$\sgn_\bH=(-1)^{\text{$\Fq$-rank of $\bH$}}$ and where we denote
by $|\bH^F|_{q'}$ the part prime to $q$ of $|\bH^F|$.
\end{enumerate}
\end{lemma}
\begin{proof}
The formula in  (i) is well known (see, \eg,  \cite[1.8]{Le}). Part (ii)
is obtained directly  from (i), taking into account  the following three
facts. First, it follows from  \cite[6.5(i)]{Sp} that the eigenvalues of
$\phi$  on  $Y(\bT)\otimes\bbR$  are the  $\delta_i\inv$;  secondly,  if
$\delta_i\neq  \delta_i\inv$, the  corresponding  basic invariants  have
the  same  degree. The  latter  fact  follows  because $\phi$  is  real,
and  so its  eigenvalues  come  in conjugate  pairs.  As a  consequence,
we have $\prod_i(q^{d_i}-\delta_i)=\prod_i(q^{d_i}-\delta_i\inv)$, which
is   required  for   the   identity  (ii).   Finally,   one  needs   the
fact  that  $\sgn_\bH=\det_{Y(\bT)\otimes\bbR}(-\phi)$   
which holds because  for  any
automorphism  $\phi$  of  finite  order   of  a  lattice  $Y$,  we  have
$\det_{Y\otimes\bbR}(\phi)=(-1)^d$ where  $d$ is the codimension  of the
fixed point subspace of $\phi$ in $Y\otimes\bbR$.
\end{proof}
\begin{remark}\label{star}
In this  work, we  shall encounter  several functions,  whose definition
generally involves  the number  of $F$-fixed points  of some  variety on
which $F$ acts,  and which are (Laurent) polynomials in  $q$. This means
that  if  $\phi$  remains  fixed  but  $q$ is  allowed  to  vary  as  in
\ref{orderinv}, they  are Laurent polynomials  in $q$. Examples  of such
functions  include the  orders of  $\Fq$-groups (as  in \ref{orderinv}),
$\tilde  P_{\ci,\ck}$, and  for a  unipotent element  $u\in \bG^F$  with
a  fixed  parametrization  (\eg,  in  the  Bala-Carter  classification),
$Q_{wF}(u)$, and $|C_\GF(u)|$.  In the case of functions  in $\CuG$, the
term polynomial  will be used when  they are linear combinations  of the
$\CY_\ci$, with coefficients  which are polynomials in  the above sense.
For any  such function $f(q)$, we  use the notation $f^*$  to denote the
function defined by  $f^*(q)=f(q\inv)$. The $\CY_\ci$ are  fixed by this
operation.
\end{remark}

The next result gives some properties of the function $\ZL\in \CWG$.

\begin{lemma}\label{zfunction}
\begin{enumerate}
\item We have
$|Z^{0vF}_\bL|=|Z^{0F}_\bG|\sum_{i=0}^l(\tilde r)^{\wedge i}(vF)q^{l-i}(-1)^i$
where  $l=\dim  Z^0_\bL-\dim  Z^0_\bG$  and  where  $\tilde  r$  is  the
restriction  to  $\WGLF$  of  the character  of  the  representation  of
$\WGL\semi\genby F$ on $Y(Z^0_\bL/Z^0_\bG)\otimes \bbR$, which is an extension
of the reflection character $r$ of $\WGL$.
\item We have $\ZL(q\inv)=\sgn_{Z_ \bL}q^{-\dim Z_\bL}\te\cdot\ZL(q)$.
\end{enumerate}
\end{lemma}
\begin{proof}
We have
$|Z^{0vF}_\bL|=|Z^{0F}_\bG|
\sum_i(-1)^i\Trace(vF|H^i_c(Z^0_\bL/Z^0_\bG))$. As in \cite[proof of 5.7]{DL}
or \cite[(1.4)]{Le}, we have
$$|Z^{0vF}_\bL|=|Z^{0F}_\bG|\sum_i(-1)^i
q^{l-i}\Trace(vF,\wedge^iY(Z^0_\bL/Z^0_\bG))$$
where $l=\dim Z^0_\bL-\dim Z^0_\bG$.
Now the  space $Y(Z^0_\bL/Z^0_\bG)\otimes \bbR$ realizes  the reflection
representation  of  the  Coxeter  group  $\WGL$, as  can  be  seen  from
\cite[9.2]{ICC}  and \cite[theorem  6]{Ho}, and  part (i)  of the  lemma
follows.

For  (ii),  let $v\in  \WGL$  and  consider  the torus  $Z_\bL^0$,  with
Frobenius action  $vF$. From \ref{orderinv}  (ii) applied here,  we have
$|Z^{0vF}_\bL|(q\inv)=\sgn'_{Z_\bL}q^{-\dim  Z_\bL}  |Z^{0vF}_\bL|(q)$,
where  $\sgn'_{Z_\bL}$  is  the  $\Fq$-rank of  $Z_\bL$  with  Frobenius
$vF$.  But, since  $\te(vF)=\det_{Y(Z^0_\bL)}(v)$ (recall that $v$  acts
trivially on $Z_\bG$  and that $\te$ is the trivial  extension), we have
$\sgn'_{Z_\bL}=\sgn_{Z_\bL} \te(vF)$.
\end{proof}

When  $\bG$ is  quasi-simple,  $\WGL$  is irreducible,  so  that $r$  is
irreducible. We then have
\begin{lemma}\label{preferred}
When  $r$ is irreducible, $\tilde  r$ is the preferred  extension of the
reflection character.
\end{lemma}
\begin{proof}
The  lemma is a consequence of  the definition  of  the  preferred  extension
in   \cite[17.2]{CS},  and   the   fact  (which   can   be  checked   by
tracing  through   \cite[9.2]{ICC})  that  if  we   write  $F=q\phi$  on
$V=Y(Z^0_\bL/Z^0_\bG)\otimes\bbR$ so  that $\tilde  r$ is  the extension
of $r$ in which $F$ acts  via $\phi$, the automorphism $\phi$ stabilizes  a set of
positive roots of a root system for $\WGL$ in $V$. We need just consider
the case  when $\phi$ is non-trivial,  so that $(\WGL,\phi)$ is  of type
$\lexp 2A_n$,  $\lexp2E_6$, $\lexp  3D_4$ or  $\lexp2D_n$. In  the cases
$\lexp 2A_n$,  $\lexp2E_6$, in the  language of \cite[17.2]{CS}  one has
$a_r=1$ so the  preferred extension is the one where  $F$ acts by $-w_0$,
which  agrees with  $\phi$.  In  the case  $\lexp  3D_4$, the  preferred
extension  is  the  only  rational  one so  again  agrees  with  $\phi$.
Finally,  in the  case $\lexp2D_n$  one checks  from the  description in
\cite[17.2]{CS} that the  preferred extension is the  one which realizes
the  reflection  representation  of $B_n\simeq  D_n\semi\genby  F$,  and
indeed $\phi$ acts  as a reflection, since it acts by  exchanging two of the
simple roots and fixing the others.
\end{proof}
If $\bG$ is not quasi-simple the group $\WGL$ is a direct product of the
irreducible Coxeter groups $W_{\bG_i}(\bL)$  where $\bG_i$ runs over the
quasi-simple  components  of  $\bG$.  The representation  of  $\WGL$  on
$Y(Z^0_\bL/Z^0_\bG)\otimes  \bbR$  decomposes  into  the  sum  over  $i$
of  summands  isomorphic  to  the  reflection  representation  $r_i$  of
the  component  $W_{\bG_i}(\bL)$ on  $Y(Z^0_\bL/Z^0_{\bG_i})\otimes\bbR$
tensored with the identity representations  of the other components. The
action  of $F$  permutes  the $r_i$  in  the same  way  it permutes  the
$\bG_i$. Since the preferred extension  of the identity is the identity,
it follows that  if $\bG_i$ is $F$-stable, the extension  of $r_i$ which
appears in  $Y(Z^0_\bL/Z^0_\bG)\otimes \bbR$ is the  preferred extension
of $r_i$.

\smallskip
We now  describe Lusztig restriction  in terms of the  generalized Green
functions, which form  a basis of the space $\CuG$.  Let $w\in \WGL$ and
suppose  $\bM$ is  a rational  Levi subgroup  which contains  a rational
conjugate  $\bL_w$  of $\bL$.  Then  we  shall use  the  identifications
explained  in \ref{twist}, \ref{split}  to consider  $Q^\bM$ as  a linear
isomorphism between $\CC(\WMLF)$ and $\CC_\CIM(\MF)$.

\begin{theorem}\label{indres} Let  $\bM$ be a rational  Levi subgroup of
some parabolic  subgroup of $\bG$. Then  $\sR^\bG_\bM\circ Q^\bG=0$
unless $\bM$  contains some  rational $\bG$-conjugate $\bL_w$  of $\bL$,
and if this condition holds, then in the above notation, we have
\begin{enumerate}
\item $\sR^\bG_\bM\circ Q^\bG=Q^\bM\circ\Res^\WGLF_\WMLF$.
\item $R^\bG_\bM\circ Q^\bM=Q^\bG\circ\Ind^\WGLF_\WMLF$.
\end{enumerate}
\end{theorem}

\begin{proof}
We need  only verify  the statements  on a basis  of the  relevant space
of  functions. We  start  by  proving (ii),  for  which  it suffices  to
evaluate both sides on $\tX_\ci$ for $\ci\in \CIM^F$. By Frobenius
reciprocity, \ref{a} (i) can be written as
$$\displaylines{R_\bM^\bG(\tX_\ci)=\sum_{\ck\in\CI^F}
\scal{\Ind^\WGLF_\WMLF\tphi_\ci}{\tphi_\ck}\WMLF Q^\bG(\tphi_\ck)=\hfill\cr
\hfill Q^\bG(\sum_{\ck\in\CI^F}\scal{\Ind^\WGLF_\WMLF\tphi_\ci}{\tphi_\ck}
\WMLF\tphi_\ck)=Q^\bG(\Ind^\WGLF_\WMLF\tphi_\ci),}$$
whence (ii) follows.

Now   take   $\theta\in\CWG$   and   consider   $\sR^\bG_\bM\circ
Q^\bG(\theta)$. The space $\Cuni(\MF)$  has a basis
$\mathop{\bigcup}\limits_{\CI'_\bM} \{\tX_\ci\mid
\ci\in {\CI'_\bM}^F\}$ where $\CI'_\bM$  runs over the $F$-stable blocks
of $\bM$. Now
$$\scal{\sR^\bG_\bM\circ Q^\bG(\theta)}{\tX_\ci}\MF=
\scal{Q^\bG(\theta)}{R_\bM^\bG(\tX_\ci)}\GF,$$
and   by  \ref{a}   the   function   $R_\bM^\bG(\tX_\ci)$  is   in
$\CC_{\CI'_\bG}(\GF)$,   where  $\CI'_\bG$   is  the   block  of   $\bG$
corresponding  to  $\CI'_\bM$.  Thus  the   scalar  product  is  $0$  if
$\CI'_\bG$ is not equal to $\CI$. Furthermore, the block $\CI$ is of the
form  $\CI'_\bG$  for some (unique by \cite[1.2]{DLM2})
block  $\CI'_\bM$  of  $\bM$ only  if  $\bM$
contains a $\bG$-conjugate $\bL_w$ of  $\bL$, whence the first statement
of the theorem.

It follows  also, that  to prove  (i), we  need only  show that  for any
$\theta\in\CWG$,  if  we  apply  both  sides of  (i)  to  $\theta$,  the
resulting  functions  have the  same  inner  product with  any  function
in  $\CC_\CIM(\MF)$. But  $\CC_\CIM(\MF)$  is spanned  by the  functions
$Q^\bM(\psi)$ with $\psi\in\CC(W_\bM(\bL_w).F)$, so  that it suffices to
consider inner products with these functions. We have

\begin{align*}
\scal{\sR^\bG_\bM\circ Q^\bG(\theta)}{Q^\bM(\psi)}\MF
&=\scal{Q^\bG(\theta)}{R^\bG_\bM(Q^\bM(\psi))}\GF&&\\
&=\scal{Q^\bG(\theta)}{Q^\bG\circ\Ind^\WGLF_\WMLF(\psi)}\GF&&\text{ by (ii)}\\
&=\scal{\theta \ZL\inv}{\Ind^\WGLF_\WMLF(\psi)}\WGLF
&&\text{by \ref{scalar product QG}}\\
&=\scal{\ZL\inv \Res^\WGLF_\WMLF(\theta)}\psi\WMLF&&\\
&=\scal{Q^\bM\circ\Res^\WGLF_\WMLF(\theta)}{Q^\bM(\psi)}\MF,
\end{align*}
which completes the proof.
\end{proof}

\begin{remark}\label{isom} The result \ref{indres} may be expressed as
asserting the  commutativity of the following diagrams.
$$\begin{CD}
\CWG@>{Q^\bG}>>\CJG\\@A\Ind AA@A{R^\bG_\bM}AA\\
\CC(\WMLF)@>{Q^\Mo}>>\CC_{\CI_\Mo}(\Mo^{\dot wF})
\end{CD}$$
and
$$\begin{CD}
\CWG@>{Q^\bG}>>\CJG\\@V\Res VV@V{\sR^\bG_\bM}VV\\
\CC(\WMLF)@>{Q^\Mo}>>\CC_{\CI_\Mo}(\Mo^{\dot wF})
\end{CD}$$
\end{remark}

As an immediate corollary, we have the following explicit formula for the
Lusztig restriction of the generalized Green functions.

\begin{corollary}\label{*RMG(Q)} With notation as in \ref{indres}, we have
$$\sR^\bG_\bM Q^\bG_{vF}=|\WML|\inv\sum_{\{x\in \WGL\mid
x(vF)x\inv\in \WMLF\}}Q^\bM_{x(vF)x\inv}.$$
\end{corollary}
\begin{proof} It is easy to see that
$$ \Res^\WGLF_\WMLF\gamma_{vF}=|\WML|\inv
\sum_{x\in \WGL, x(vF)x\inv\in \WMLF}\gamma_{x(vF)x\inv}. $$
The result now follows immediately by applying \ref{indres} (i) to
the function $\gamma_{vF}$.
\end{proof}
The duality involution  $\cD_\bG$ (restricted to $\CJG$)  has an elegant
description in this setting.
\begin{proposition}\label{DG} (\cf\ \cite{US})
Let $\cD_\bG$ be the duality involution; then
\begin{enumerate}
\item We have $\cD_\bG(Q_{wF})=\eta_\bL\te(wF)Q_{wF}$,
where, for any reductive group $\bG$ we write
$\eta_\bG=(-1)^{\text{semisimple $\Fq$-rank of \bG}}=\sgn_\bG\sgn_{Z_\bG}$.
\item    The   duality    involution   $\cD_\bG:\CJG\longrightarrow\CJG$
corresponds under $Q^\bG$ to multiplication by $\eta_\bL\te$ in $\CWG$.
In  particular $\cD_\bG(\tX_\ci)=\eta_\bL\sgn_\ci\tX_{\hat\ci}$,  where
$\hat\ci$ and $\sgn_\ci$ are defined in \ref{epsilon}.
\end{enumerate}
\end{proposition}
\begin{proof}
The statement (i) may be found in \cite[\S 8]{US} whose proof applies to the
twisted case without change. The first statement in (ii) follows immediately
since $Q^\bG$ is linear, and the second statement follows from the relation
$\te\otimes\tphi_\ci=\sgn_\ci\tphi_{\hat\ci}$ (\ref{epsilon}). 
\end{proof}

\section{Unipotently supported class functions
and twisted class functions on reflection groups}
For $\ci\in\CI^F$ define a function $\tQ_\ci$ on $\WGLF$ by

\begin{equation}\label{defQi}
\tQ_\ci(wF)=\frac1{a_\ci}\sum_{a\in A(u)}
q^{-c_\ci}\overline{\CY_\ci(u_a)}Q_{wF}(u_a)
\end{equation}
where  we   fix  $u\in  C_\ci^F$  and   set  $A(u)=C_\bG(u)/C_\bG^0(u)$,
$a_\ci=|A(u)|$  and   take  $u_a$   to  be   a  representative   of  the
$\GF$-orbit   in   $C^F$  which   corresponds   to   the  $F$-class   of
$a\in  A(u)$.   Then  using  the  relation
\begin{equation}\label{orth.Y}
a_\ci\inv\sum_{a\in  A(u)}
\CY_\ci(u_a)\overline{\CY_\cj(u_a)}=\delta_{\ci,\cj}
\end{equation}
(see \cite[1.5]{DLM2}) and \ref{QwF}, we obtain
\begin{equation}\label{QiwF}
\tQ_\ci(wF)=\sum_{\cj\in\CI^F}\tphi_\cj(wF)\tilde P_{\ci,\cj}
\end{equation}

Observe  that  the  definition   of  $\tQ_\ci$  given  in  (\ref{defQi})
above,  makes  sense  even  when   $\ci\not\in  \CI$,  but  then,  since
$\tilde P_{\ci,\cj}=0$  when $\ci$  and $\cj$  are in  different blocks,
$\tQ_\ci=0$.

\begin{proposition}\label{invQi} We have
$Q_{wF}(u)=\sum_{\ci\in\CI^F}\tQ_\ci(wF)\tilde\CY_\ci(u)$
\end{proposition}
\begin{proof}
As remarked above, if  $\ci\notin  \CI$ then \ref{QiwF}  shows that
the corresponding summand of the right hand side is $0$, since then
$\tilde  P_{\ci,\cj}=0$ for all $\cj\in\CI^F$. So
$$\sum_{\ci\in\CI^F}\tQ_\ci(wF)\tilde\CY_\ci(u)=
\sum_{\ci\in \CP^F}\tQ_\ci(wF)\tilde\CY_\ci(u)$$
We now use the second orthogonality formula for the $\CY_\ci(u)$:
\begin{equation}\label{2orthogY}
\sum_{\ci\in \CP^F}\CY_\ci(u)\overline{\CY_\ci(u')}=\begin{cases}
|A(u)^F|&\text{ if $u\sim_\GF u'$}\\0&\text{otherwise}\end{cases}
\end{equation}
where $\sim_\GF$ means $\GF$-conjugacy. Thus
$$\displaylines{\sum_{\ci\in \CI^F}\tQ_\ci(wF)\tilde\CY_\ci(u)=
\sum_{\ci\in \CP^F, a\in A(u)}a_\ci\inv\overline{\CY_\ci(u_a)}Q_{wF}(u_a)
\CY_\ci(u)=\hfill\cr\hfill |A(u)|\inv|A(u)^F|
\#\{a\mid u_a\sim_\GF u\}Q_{wF}(u)=Q_{wF}(u)\cr}$$
\end{proof}

Note that  the equation \ref{2orthogY}  will often be used  when $u=u_a$
for some $a\in  A(u)$, in which case  we have $|A(u)^F|=|C_{A(u)}(aF)|$.
The functions $\tilde \CY_\ci$ form a basis of $\CJG$ as $\ci$ runs over
$\CI^F$. The next result relates the $\tQ_\ci$ to expansions in terms of
this basis.

\begin{lemma}\label{coefQ}
\begin{enumerate}
\item For any function $f\in \CJG$, the coefficient of $f$ in the basis
$\tilde \CY_\ci$ is
\begin{equation}\label{coefY}
\frac1{a_\ci}\sum_{a\in A(u)}q^{-c_\ci}\overline{\CY_\ci(u_a)}f(u_a).
\end{equation}
\item For any function $\theta\in \CWG$, we have
$Q^\bG(\theta)=\sum_{\ci\in \CI^F}\scal\theta{\tQ_\ci}\WGLF\tilde\CY_\ci$.
\item  The  functions  $(Q^\bG)\inv(\tilde\CY_\ci)$ form  the  basis  of
$\CWG$ which is dual to the basis $\{\tQ_\ci\}$.
\end{enumerate}
\end{lemma}
\begin{proof}
By  \ref{invQi}  and the  definition  (\ref{defQi})  of $\tQ_\ci$,
(i)  holds when  $f=Q_{wF}$,  and since  the $Q_{wF}$  form  a basis  of
$\CJG$  and the  formula  \ref{coefY} is  linear in  $f$,  (i) holds  in
general.  Similarly,  (ii)  holds when  $\theta=\gamma_{wF}$,  again  by
\ref{invQi}.  By linearity,  (ii) holds  generally. The  statement (iii)
follows immediately from (ii).
\end{proof}

\section{Lusztig's algorithm and orthogonality relations for generalized
Green functions}
We shall require
\begin{lemma}\label{invscal}
Let  $H$  be  a  finite group,  $\chi_1,\chi_2,\cdots$  the  irreducible
characters of  $H$ (over  a field  of characteristic  zero) and  $f$ any
class function  on $H$  which is  non-zero at each  element of  $H$. Let
$f\inv$  be the  pointwise  inverse  of $f$.  Then  we  have the  matrix
equation
\begin{equation}\label{scalinv}
\{\scal{f\inv\chi_i}{\chi_j}H\}_{i,j}=\{\scal{f\chi_i}{\chi_j}H\}_{i,j}\inv.
\end{equation}
\end{lemma}
\begin{proof} Since the $\chi_i$ form  an orthonormal basis of the space
of class functions on $H$, the  left side of \ref{scalinv} is simply the
matrix  of  the  linear  transformation  induced  by  multiplication  by
$f\inv$,  and  the  assertion  is  no more  than  the  observation  that
multiplication by $f\inv$ is the inverse of multiplication by $f$.
\end{proof}
The statement  \ref{invscal} remains valid  when $H$ is a  finite coset,
the $\chi_i$ are extensions to $H$  of the irreducible characters of the
underlying group, and $f$ is a twisted class function on $H$.

We now recall the algorithm outlined  by Lusztig in \cite[\S 24]{CS} for
the computation  of the  polynomials $P_{\ci,\ck}$:  Lusztig's algorithm
is  based  on the  following  matrix  equation,  which is  an  immediate
consequence of the relation $\tX_\ci=\sum_\ck
\tilde P_{\ck,\ci}\tilde \CY_\ck$ and \ref{scalX}.
$$\lexp t{\tilde P}\tilde\Lambda\tilde P=
\{\scal{\tX_\ci}{\tX_\ck}\GF\}_{\ci,\ck}
=\{\scal{\ZL\inv \tphi_\ci}{\tphi_\ck}\WGLF\}_{\ci,\ck}$$
where $\tilde P=\{\tilde P_{\ci,\ck}\}_{\ci,\ck}$
and $\tilde\Lambda=\{\scal{\tilde \CY_\ci}{\tilde \CY_\ck}\GF\}_{\ci,\ck}$.
We shall use the inverse of this equation:
$${\tilde P}\inv \tilde\Lambda\inv(\lexp t{\tilde P}\inv)=\tilde\Omega$$
where $\tilde\Omega=\{\tilde\omega_{\ci,\ck}\}_{\ci,\ck}$ and
$\tilde\omega_{\ci,\ck}=
\scal{\ZL\tphi_\ci}{\tphi_\ck}\WGLF$,
the inverse of  the matrix on the  right hand side being  given by Lemma
\ref{invscal}. The  matrix $\tilde\Omega$  may be considered  known (see
\ref{defCZGL}) since it  is given in terms of Weyl  group data. The rows
and columns  of $\tilde\Lambda$ and $\tilde  P$ may be ordered  in a way
compatible  with the  order on  unipotent classes;  they may  further be
ordered so that pairs with the same support form a connected sequence in
the  order.  Then  $\tilde\Lambda$  is  block-diagonal  and  $\tilde  P$
block-triangular with identity diagonal blocks, the blocks corresponding
to unipotent  classes. Given  $\tilde\Omega$, there are  unique matrices
$\tilde\Lambda$ and  $\tilde P$  of this shape  which satisfy  the above
equation.

We note for future reference that \ref{zfunction} immediately gives
\begin{equation}\label{tildeomegaik}
\tilde\omega_{\ci,\ck}=|Z^{0F}_\bG|
\sum_{i=0}^lq^{l-i}(-1)^i\scal{\tphi_\ci\otimes\tphi_\ck}
{{\tilde r}^{\wedge i}}\WGLF
\end{equation}
where  $l=\dim  Z^0_\bL-\dim  Z^0_\bG$  and  where  $\tilde  r$  is  the
restriction  to  $\WGLF$  of  the character  of  the  representation  of
$\WGL\semi\genby  F$ on  $Y(Z^0_\bL/Z^0_\bG)\otimes \bbR$,  which is  an
extension of the reflection character $r$ of $\WGL$.

The following proposition is a generalization of \cite[1.1.4]{Ka}.

\begin{corollary}\label{orthQi}
(second orthogonality formula for Green functions)
$$\scal{\ZL\tQ_\ci}{\tQ_\cj}\WGLF=
\{\scal{\tilde\CY_\ci}{\tilde\CY_\cj}\GF\}\inv_{\ci,\cj}
=\begin{cases}
a_\ci\inv \sum\limits_{a\in A(u)}\dfrac{|C^0_\bG(u_a)^F|}
{q^{2c_\ci}}\CY_\ci(u_a)\overline {\CY_\cj(u_a)}
&\text{if $C_\ci= C_\cj$}\\0&\text{otherwise}
\end{cases}$$
where notation is as in \ref{defQi}.
\end{corollary}
\begin{proof} Using  the values  given in \ref{QiwF}  for $\tQ_\ci$
and $\tQ_\cj$, we obtain:
$$\{\scal{\ZL\tQ_\ci}{\tQ_\cj}\WGLF\}_{\ci,\cj}=
\tilde P\{\scal{\ZL\tphi_\ck}{\tphi_{\ck'}}\WGLF\}_{\ck,\ck'}
\lexp t{\tilde P}=\tilde P\tilde\Omega\lexp t{\tilde P}=
\{\scal{\tilde\CY_\ci}{\tilde\CY_\cj}\GF\}\inv_{\ci,\cj}
$$
Now $\scal{\tilde\CY_\ci}{\tilde\CY_\cj}\GF$ is $0$ if $C_\ci\ne C_\cj$
and otherwise is equal to
\begin{equation}\label{scalY}
\sum_{a\in H^1(F,A(u))} |C_\GF(u_a)|\inv\tilde\CY_\ci(u_a)
\overline{\tilde\CY_\cj(u_a)}=a_\ci\inv\sum_{a\in A(u)}
|C^0_\bG(u_a)^F|\inv\tilde\CY_\ci(u_a)\overline{\tilde\CY_\cj(u_a)}.
\end{equation}
To see \ref{scalY}, note that $(A(u_a),F)$ is isomorphic to $(A(u),aF)$,
so  that  $$|C_\GF(u_a)|=|C_{A(u)}(aF)||C^0_\bG(u_a)^F|.$$  Finally,  it
follows  from  \ref{orth.Y}  and  \ref{invscal} that  the  matrix  whose
$(\ci,\cj)$ entry  is either side of  \ref{scalY} is the inverse  of the
matrix whose $(\ci,\cj)$ entry is the expression in the statement.
\end{proof}
This in turn gives an orthogonality formula for the $Q_{wF}$, regarded
as elements of $\CWG$ for a fixed value of the argument:
\begin{corollary}\label{orthQw}
For  $u$   a  unipotent   element  of   $\GF$,  define   the  function
$Q_-(u)\in\CWG$ by $Q_-(u)(wF)=Q_{wF}(u)$ (for $wF\in\WGLF$).  Then
$$\displaylines{\scal{Q_-(u)}{\ZL Q_-(u')}\WGLF=\hfill\cr\hfill
\begin{cases}
\displaystyle |A(u)|\inv\sum_{a\in A(u)}|C^0_\bG(u_a)^F|
(\sum_{\ci\in\CI^F}\overline{\CY_\ci(u_a)}\CY_\ci(u))
(\sum_{\ci\in\CI^F}\CY_\ci(u_a)\overline{\CY_\ci(u')})
&\text{if $u\sim_\bG u'$}\\
0& \text{otherwise.}
\end{cases}
\cr}$$
\end{corollary}
\begin{proof} Applying \ref{invQi} and then \ref{orthQi}
to the left-hand side we get
$$
\displaylines{\scal{Q_-(u)}{\ZL Q_-(u')}\WGLF=
\scal{\sum_\ci\tQ_\ci\tilde\CY_\ci(u)}{\sum_\cj\ZL\tQ_\cj\tilde\CY_\cj(u')}\WGLF
\hfill\cr\hfill
=\sum_{\ci,\cj}\tilde\CY_\ci(u)\overline{\tilde\CY_\cj(u')}
\{\scal{\tilde\CY_\ci}{\tilde\CY_\cj}\GF\}\inv_{\ci,\cj}.\cr}
$$
we then use that the matrix 
$\{\scal{\tilde\CY_\ci}{\tilde\CY_\cj}\GF\}_{\ci,\cj}$
is real to write the complex conjugate of the expression in \ref{orthQi}
and we get the result.
\end{proof}
If we  sum formula \ref{orthQw} over  all blocks, we obtain  the simpler
expression:
\begin{proposition}\label{allorth}
$$\sum_\CI\scal{Q^\CI_-(u)}{\ZL Q^\CI_-(u')}\WGLF=\begin{cases}
|C_\GF(u)|&\text{if $u\sim_\GF u'$}\\
0&\text{otherwise}\end{cases}$$
where $\CI$ runs over the rational blocks and where we put a superscript $\CI$
on the $Q_-$ to show which block they come from.
\end{proposition}
\begin{proof} The sum over all blocks of the right-hand side of \ref{orthQw}
is, when $u\sim_\bG u'$
$$|A(u)|\inv\sum_{a\in A(u)}|C^0_\bG(u_a)^F|
(\sum_{\ci\in\CP^F}\overline{\CY_\ci(u_a)}\CY_\ci(u))
(\sum_{\ci\in\CP^F}\CY_\ci(u_a)\overline{\CY_\ci(u')}),$$
which, using the second orthogonality formula
\ref{2orthogY} for $\CY_\ci$ reduces to
$$|A(u)|\inv\sum_{\{a\in A(u)\mid u_a\sim_\GF u\text{ and }u_a\sim_\GF u'\}}
|C^0_\bG(u_a)^F||A(u)^F||A(u')^F|$$
which is $0$ unless $u \sim_\GF u'$ and equal to $|C_ \GF(u)|$ otherwise.
\end{proof}
\section{Gelfand-Graev characters and their Lusztig restriction}
As in \cite{US}  and \cite{DLM2}, for $\ci\in\CI^F$  and $u\in C_\ci^F$,
we  define $\Gamma_\ci=\sum_{a\in  A(u)}\CY_\ci(u_a)\Gamma_{u_a}$, where
$\Gamma_{u_a}$ is  the generalized  Gelfand-Graev character  attached to
the class of $u_a$, and other  notation is as in \ref{QiwF}. 
\begin{proposition}\label{Gammatilde} 
We  have  $\Gamma_\ci=a_\ci\zeta_\CI\inv Q^\bG(\te\ZL\tQ^*_\ci)$,  where
$\zeta_\CI$ is  a fourth root of  unity (the one associated  to $\CI$ in
\cite[7.2]{US} when $\bG$ is split).
\end{proposition}
\begin{proof}
We start from  the formula \cite[7.5 (b)]{US} of Lusztig,  which must be
modified for the case  a non-split  group  in a  way hinted at in
\cite[8.7]{US}. We claim  that for a possibly  non-split group, the
equation \cite[7.5(b)]{US} should read
\begin{equation}\label{GammaLusztig}
\Gamma_{\ci_0}=a_{\ci_0}\zeta_\CI\inv\sum_{\ci,\ci_1}|\WGL|\inv\sum_{w\in\WGL}\tphi_{\hat\ci_1}(wF)
\tphi_\ci(wF)|Z_\bL^{0wF}|\tilde P^*_{\ci_0,\ci}\sgn_{\ci_1}\tX_{\ci_1}.
\end{equation}
The  only  part  of the  generalization which is not obvious, and which  is
the source of the coefficient $\sgn_{\ci_1}$  in the above formula, is
(as indicated in \cite[8.7]{US}) the lemma  \cite[7.2]{US} whose statement
should be changed for the general situation to read
$\hat\tX_\ci\mid_{\GF_{\text{uni}}}=\zeta_\CI
q^{(\dim\bG-\dim  Z_\bL)/2}\sgn_\ci\tX_{\hat\ci}$.  The proof  given  in
\cite[7.2]{US} cannot be applied in  our  more general  case, since  $\dim
V_\ci$ has to be replaced by $\Trace(F\mid V_\ci)$, which might vanish.
Nonetheless the generalization may be proved  by considering  a Frobenius  twisted by
various $v\in\WGL$ on  the induced sheaf which Lusztig  considers in that
proof.

We now rewrite \ref{GammaLusztig} as
$$\begin{aligned}
\Gamma_{\ci_0}&=a_{\ci_0}\zeta_\CI\inv|\WGL|\inv\sum_{w\in\WGL}|Z_\bL^{0wF}|
\sum_\ci\tphi_\ci(wF)\tilde P^*_{\ci_0,\ci}
\sum_{\ci_1}\tphi_{\hat\ci_1}(wF)\sgn_{\ci_1}\tX_{\ci_1}\\
&=a_{\ci_0}\zeta_\CI\inv|\WGL|\inv\sum_{w\in\WGL}|Z_\bL^{0wF}|
\sum_\ci\tphi_\ci(wF)\tilde P^*_{\ci_0,\ci}
\sum_{\ci_1}\te(wF)\tphi_{\ci_1}(wF)\tX_{\ci_1}\text{ by \ref{epsilon}}\\
&=a_{\ci_0}\zeta_\CI\inv|\WGL|\inv\sum_{w\in\WGL}|Z_\bL^{0wF}|
\sum_\ci\tphi_\ci(wF)\tilde P^*_{\ci_0,\ci}
\te(wF)Q_{wF}\text{ by \ref{QwF}}\\
&=a_{\ci_0}\zeta_\CI\inv|\WGL|\inv\sum_{w\in\WGL}|Z_\bL^{0wF}|\te(wF)
\tQ_{\ci_0}^*
Q_{wF}\text{ by \ref{QiwF}}.\\
\end{aligned}$$
The proposition now follows by \ref{defqg}(iii).
\end{proof}

Let us write $\tilde\zeta_\CI$ for the root of unity denoted by $\zeta_\CI$
in \cite{DLM2}. The point of this notation is to distinguish $\tilde\zeta_\CI$ 
and $\zeta_\CI$, since they turn out to be different generalizations to non-split groups 
of Lusztig's constant.
\begin{proposition}\label{scalgammaX}
For any reductive group, let $\sigma_\bG:=(-1)^{\text{semi-simple
rank$(\bG)$}}$. Then
$\tilde\zeta_\CI=\eta_\bL\sigma_\bL\zeta_\CI$.
\end{proposition}
\begin{proof}
We have
$$\begin{aligned}
\scal{D_\bG\Gamma_\ci}{\tX_\ck}\GF
&=\scal{\Gamma_\ci}{D_\bG\tX_\ck}\GF
=\eta_\bL\sgn_\ck\scal{\Gamma_\ci}{\tX_{\hat\ck}}\GF\text{ by  \ref{DG}(ii)}\\
&=\eta_\bL\sgn_\ck\scal{\ZL\inv(Q^\bG)\inv(\Gamma_\ci)}
{(Q^\bG)\inv(\tX_{\hat\ck})}\GF
\text{ by \ref{scalar product QG}}\\
&=\eta_\bL\sgn_\ck a_\ci\zeta_\CI\inv\scal{\te\tQ_\ci^*}{\tphi_{\hat\ck})}\WGLF\\
&=\eta_\bL a_\ci\zeta_\CI\inv\scal{\tQ_\ci^*}{\tphi_\ck)}\WGLF\text{ by \ref{epsilon}}\\
&=\eta_\bL a_\ci\zeta_\CI\inv\tilde P^*_{\ci,\ck}\text{ by the $*$ of \ref{QiwF}.}
\end{aligned}$$
The  equation \cite[1.7]{DLM2} is transformed into this last relation
if  $\sigma_\bL\tilde\zeta_\CI\inv$  
is replaced by $\eta_\bL\zeta_\CI\inv$, whence the proposition.
\end{proof}

It will be convenient to  use the normalization $\tG_\ci=a_\ci\inv\zeta_\CI\Gamma_\ci$.
We shall now discuss orthogonality relations
among the  $\tG_\ci$ and among  the $\Gamma_u$,  as well as  the Lusztig
restriction  of  the  $\tG_\ci$.  Note  that  from  \ref{Gammatilde}  it
follows  that if  $\CI$  is  a rational  block  and $\ci\in\CI^F$,  then
$\tG_\ci\in\CJG$.

\begin{lemma}\label{gammauI} For any rational block $\CI$
define $\tG_u^\CI=\sum_{\ci\in\CI}\overline{\tilde\CY_\ci^*(u)}\tG_\ci$.
If  there  is  a  pair $\ci\in\CI^F$  whose  support  contains  $u$,
the  orthogonal projection  of  $\Gamma_u$ onto  $\CJG$ is  $\zeta_\CI\inv
q^{c_\ci}\tG_u^\CI$; otherwise it is $0$.
\end{lemma}
\begin{proof}
Using \ref{2orthogY}, the defining relation for $\Gamma_\ci$ can be inverted
to give
$$
\Gamma_u=|A(u)|\inv\sum_{\ci\in\CP^F}\overline{\CY_\ci(u)}\Gamma_\ci.
$$
If we  restrict the  above sum  to $\ci\in\CI^F$  we obtain the orthogonal
projection of $\Gamma_u$ onto $\CJG$,  since the various spaces $\CJG$ are
mutually orthogonal. The  lemma now follows  in straightforward fashion 
from the definitions.
\end{proof}
\begin{proposition}\label{Gammai} We have 
$\tG_u^\CI=Q^\bG(\te \ZL \overline {Q_-}^*(u))$.
\end{proposition}
\begin{proof}
Apply $(Q^\bG)\inv$ to the expression \ref{gammauI}
for $\tG^\CI_u$ to get
$$
(Q^\bG)\inv(\tG_u^\CI)=\sum_{\ci\in\CI^F}\overline{\tilde \CY^*_\ci(u)}
\te\ZL\tQ^*_\ci.
$$
Now take the complex conjugate of the $^*$ of the relation \ref{invQi} and
substitute into  this last equation.  Taking
into account that  the functions $\tQ_\ci$ are real  valued (i.e. stable
under  complex conjugation), which is a consequence of  \ref{QiwF} since  the $\tphi_\ci$  are
real, we  obtain the proposition.
\end{proof}
\begin{corollary}\label{scalgammai}
We have $\scal{\tG_\ci}{\cD_\bG\tG_\ck}\GF=\sgn_\bG q^{\dim Z_\bL}
\left(\{\scal{\tilde\CY_\ci} {\tilde\CY_\ck} \GF\}\inv_{\ci,\ck}\right)^*$,
which is zero if $C_\ci\neq C_\ck$. 
\end{corollary}
\begin{proof}
We have 
$$\begin{aligned}
\scal{\tG_\ci}{\cD_\bG\tG_\ck}\GF&=
\scal{\ZL\inv(Q^\bG)\inv(\tG_\ci)}{\eta_\bL\te(Q^\bG)\inv(\tG_\ck)}\WGLF
\text{ by \ref{scalar product QG} and \ref{DG}(ii)}\\
&=\eta_\bL\scal{\tQ^*_\ci}{\ZL\te\tQ^*_\ck}\WGLF\\
&=\eta_\bL\sgn_{Z_\bL}q^{\dim Z_\bL}\scal{\tQ^*_\ci}{\ZL^*\tQ^*_\ck}\WGLF
\text{ by \ref{zfunction}}\\
&=\sgn_\bL q^{\dim Z_\bL}\scal{\tQ^*_\ci}{\ZL^*\tQ^*_\ck}\WGLF
\text{ since $\eta_\bL=\sgn_\bL\sgn_{Z_\bL}$}\\
&=\sgn_\bL q^{\dim Z_\bL}\left(\scal{\tQ_\ci}{\ZL\tQ_\ck}\WGLF\right)^*\\
&=\sgn_\bL q^{\dim Z_\bL}
\left(\{\scal{\tilde\CY_\ci}{\tilde\CY_\ck}\GF\}\inv_{\ci,\ck}\right)^*
\text{ by \ref{orthQi}}.\\
\end{aligned}
$$
The result now follows because $\sgn_\bL=\sgn_\bG$ since $\bL$ is $\bG$-split.
\end{proof}

\begin{corollary}\label{scalgammauI}
Let $u,v\in \GF$ be unipotent elements and $\CI$ a rational block. Then  
$$\scal{\tG_u^\CI}{\cD_\bG\tG_v^\CI}\GF=\sgn_\bG q^{\dim Z_\bL}
\left(\scal{\overline Q_-(u)}{\ZL\overline  Q_-(v)} \WGLF\right)^*,$$
which is non-zero only if $u$ and $v$ are conjugate in $\bG$.
\end{corollary}
\begin{proof}
We have, from  \ref{Gammai}, proceeding as in \ref{scalgammai}
$$\begin{aligned}
\scal{\tG_u^\CI}{\cD_\bG\tG_v^\CI}\GF
&=\scal{\te\overline Q_-^*(u)}{\eta_\bL\te^2\ZL\overline Q_-^*(v)}\WGLF\\
&=\eta_\bL\scal{\te\overline Q_-^*(u)}{\ZL\overline Q_-^*(v)}\WGLF\\
&=\eta_\bL\sgn_{Z_\bL} q^{\dim Z_\bL}
\scal{\te\overline Q_-^*(u)}{\te \ZL^*\overline Q_-^*(v)}\WGLF\\
&=\sgn_\bL q^{\dim Z_\bL}\scal{\overline Q_-^*(u)}{\ZL^*\overline Q_-^*(v)}\WGLF\\
&=\sgn_\bL q^{\dim Z_\bL}\left(\scal{\overline Q_-(u)}{\ZL\overline Q_-(v)}\WGLF\right)^*,\\
\end{aligned}
$$
and the  result follows  as in  \ref{scalgammai}. The  last remark  is a
consequence of the evaluation of the right side in \ref{orthQw}.
\end{proof}
\begin{corollary}\label{orthgammau}
For any  pair $u,v$ of  unipotent elements of  $\GF$, we have
$$\scal{\Gamma_u}{D_\bG\Gamma_v}=
\begin{cases}
\sgn_\bG \sgn_{C_\bG(u)}|C_\GF(u)|_{q'}&\text{ if $u\sim_\GF v$}\\
0&\text{ otherwise.}\\
\end{cases}
$$
\end{corollary}
\begin{proof}  From \ref{gammauI}, we see 
$\scal{\Gamma_u}{D_\bG\Gamma_v}=
\sum_\CI\scal{\zeta_\CI\inv q^{c_\ci}\tG_u^\CI}
{\zeta_\CI\inv q^{c_\cj}\cD_\bG\tG_v^\CI}\GF$
where the  sum is over all blocks which contain two pairs $\ci$, $\cj$
whose support respectively contains $u$ and $v$.
By \ref{scalgammauI} this sum is $0$ if $u$ and $v$ are not $\bG$-conjugate;
otherwise we obtain
$$\scal{\Gamma_u}{D_\bG\Gamma_v}=
\sgn_\bG q^{\codim(\class(u))}\sum_\CI
\left(\scal{\overline Q^\CI_-(u)}{\ZL\overline Q^\CI_-(v)}\WGLF\right)^*.$$
We now apply \ref{allorth} and \ref{orderinv}(ii) to complete the proof.
\end{proof}

\smallskip
To describe the Lusztig restrictions of the $\tG_\ci$, we  shall require
the  notion  of  ``$\Fq$-rank  relative  to  a  block'',  which  we  now
define.  Suppose $\bM$  is a  rational Levi  subgroup as  in \ref{twist}
and  \ref{split}.  If $\bT_0$  is  a  maximally split  rational  maximal
torus  of $\bL$  (and hence  of  $\bG$), the  coset $W_\bL(\bT_0).w$  of
the  Weyl group  $W_\bL(\bT_0)$ is  uniquely  defined by  $\bM$ and  the
conditions on $\bL_w$ up to  $F$-conjugacy in $W_\bG(\bT_0)$ (see, \eg,
\cite[4.3]{DM2}). In the coset $W_\bL(\bT_0).w$ the elements whose fixed
points on  $Y(\bT_0)\otimes \bbR$ have  maximal dimension form  a single
class under $F$-conjugacy by $W_\bG(\bT_0)$, and it is the case that $w$
is among these  elements. The $\Fq$-rank of $\bM$ relative  to the block
$\CI$ is defined as the dimension of the subspace of $w$-fixed points of
$Y(Z_\bL^0)\otimes \bbR$.

\begin{definition}\label{blockrank}  
With    notation    as    in     the    previous    paragraph,    define
$$\sgn_\CI(\bM):=\sgn^\bG(w).$$ It follows from  the remarks in the last
paragraph that the right side depends only on (the $\GF$-conjugacy class
of) $\bM$.
\end{definition}
\begin{lemma}\label{restrictionsRsgn}
\begin{enumerate}
\item In the notation of Remarks \ref{twist} and \ref{split},
there exist Laurent polynomials $R_{\ci,\cj}$ in $q$
($\ci\in\CI^F$ and $\cj\in\CIM^F$) such that
$\Res^\WGLF_\WMLF\tQ_\ci=\sum\limits_{\cj\in\CIM^F}R_{\ci,\cj}\tQ_\cj$.
We have $R_{\ci,\cj}=0$ unless
$\overline{C_\cj}\subset\overline{C_\ci}\subset\overline{\Ind_\bM^\bG C_\cj}$.
\item Maintaining the above notation, we have
$\Res^\WGLF_\WMLF\te^\bG=\sgn_\CI(\bM)\te^\bM$, where $\sgn_\CI(\bM)$
is defined in \ref{blockrank}.
\end{enumerate}
\end{lemma}
\begin{proof}
 Let  $R$  be  the  matrix  with  $(\ci,\cj)$  coefficient
$R_{\ci,\cj}$ as in (i) of the statement. From \ref{QiwF}, we obtain the
matrix equation
$$\tilde P^\bG\{\scal{\tphi_\cj}{\Res^\WGLF_\WMLF\tphi_{\ci'}}\WMLF\}_{\ci',\cj}
=R\tilde P^\bM.$$
The first  statement in (i) is  now immediate, since the  entries of the
unitriangular matrix $\tilde P^\bM$  are Laurent polynomials, whence the
same is  true of its inverse.  The second statement in  (i) follows from
\ref{a} (ii).

For (ii), let $v.wF\in\WMLF$. Then 
$$ \te^\bG(v.wF)=\sgn^\bG(vw) =\sgn^\bG(v)\sgn^\bG(w)
=\te^\bM(v.wF)\sgn_\CI(\bM).  $$
\end{proof}

\begin{proposition}\label{*RMG(Gamma)}
We have
$\sR^\bG_\bM(\tG_\ci)=
\sgn_\CI(\bM)\sum_{\cj\in\CIM^F} R^*_{\ci,\cj}\tG_\cj$,
where notation is as in \ref{Gammai} and \ref{restrictionsRsgn} above.
\end{proposition}
\begin{proof}
Since $\tG_\ci\in \CJG$, it  follows from \ref{indres} that $\sR^\bG_\bM
\tG_\ci$ is  zero unless  $\bM$ contains  a rational  $\bG$-conjugate of
$\bL$. We therefore take $\bM$  as in \ref{indres}. Now by \ref{Gammai},
$\tG_\ci=Q^\bG(\te\ZL\tQ_\ci^*)$,  and  $\ZL\in   \CWG$  is  defined  in
\ref{defCZGL}. By \ref{indres}  we need only compute  the restriction to
$\WML$ of $\te  \ZL \tQ_\ci^*$, and a  straightforward calculation using
\ref{restrictionsRsgn} yields the statement.
\end{proof}

\begin{remark}\label{zeta1}
It is a consequence  of \cite[\S 2]{DLM2}
that  for  regular  blocks,  $\tilde\zeta_\CI$  is  independent  of  the
ambient  group and  the rational  structure,  \ie\ depends  only on  the
geometric data  in the cuspidal  system $(\bL,\ci_0)$. This  is asserted
without justification in  the proof of \cite[3.4]{DLM2} but  can be seen
as  follows. From  \cite[2.1]{DLM2}  and \cite[2.5]{DLM2}  one has  that
$\tilde\zeta_\CI$  is equal  (in the  notation  of {\it  loc.\ cit.})  to
$\eta_\bL\sigma_\bL\sigma^\bL_\zeta$  up to  a power  of $q$.  Using the
Hasse-Davenport relation, one  may compare the product of  Gauss sums in
\cite[2.4]{DLM2} which applies to the case of twisted  $\bL$, to  that occurring
in  a  split  group. One  finds  that  the  products  also differ  by  a
factor  $\eta_\bL\sigma_\bL$. Thus $\tilde\zeta_\CI=\tilde\zeta_\CIM$
in this case.  In particular, this applies generally to the  principal block (when
$\bL$  is a  maximal  torus). In  general, the  question  as to  whether
$\tilde\zeta_\CI=\tilde\zeta_\CIM$ in all cases  amounts to the question
of  whether  $\tilde\zeta_{\CI_\bL}$  is independent  of  the  Frobenius
structure on the triple  $(L,C_{\ci_0},\ci_0)$. Although this point does
not affect the formulation of  \ref{*RMG(Gamma)}, it is relevant to some
of the computations later in this work.
\end{remark}
\begin{remark}\label{zeta2}
The equation \ref{*RMG(Gamma)} may be expressed as follows.
$$
\sR^\bG_\bM(a_\ci\Gamma_\ci)=
\sgn_\CI(\bM)\zeta_\CI\zeta_{\CI_\bM}\inv\sum_{\cj\in\CIM^F} 
R^*_{\ci,\cj}a_\cj\Gamma_\cj
=\varepsilon_\bG\varepsilon_\bM\tilde\zeta_\CI
\tilde\zeta_{\CI_\bM}\inv\sum_{\cj\in\CIM^F} 
R^*_{\ci,\cj}a_\cj\Gamma_\cj,\\
$$
and the previous remark implies that in the regular case, the 
factor $\tilde\zeta_\CI\tilde\zeta_{\CI_\bM}\inv$ is equal to $1$.
\end{remark}
\section{Application to the regular and subregular cases}
Our  objective now  is to  apply Proposition  \ref{*RMG(Gamma)} to  some
specific cases. The general strategy will be first to compute \ref{QiwF}
explicitly in  $\bG$ and in  $\bM$ by computing certain  required values
$\tilde P_{\ci,\ck}$, and then to  use specific knowledge of restriction
of characters from $\WGLF$ to $\WMLF$.

As an example,  consider first the case when  $\ci=\creg$, where $\creg$
is a  pair in the block  $\CI$ with support the  regular unipotent class
(such a pair  is then the unique  one with regular support  in the block
$\CI$, see \cite[1.10]{DLM2}). Then the  only non-zero term in the right
hand side of formula \ref{QiwF} is $\tphi_\creg(wF)$, as
$\tilde P_{\creg,\creg}=1$ and $\tilde P_{\creg,\cj}=0$ if
$C_\creg\not\subset\overline C_\cj$. Moreover, as $\creg$
has regular support we have $\tphi_\creg=\Id$. So we get
$\tQ_\cregG=\Id_\WGLF$, whence
$\Res^\WGLF_\WMLF\tQ_\cregG=\tQ_\cregM$.
Applying \ref{zeta2} we get
$$\sR^\bG_\bM\Gamma_\cregG=\frac{a_\cregG}{a_\cregM}
\sgn_\bG\sgn_\bM\Gamma_\cregM.$$
Thus we recover lemma 3.6 of \cite{DLM2}.

\begin{proposition}\label{Psr}  Consider an  $F$-stable pair  $\cs$ with
support a  subregular class  $C_\cs$ of  $\bG$ and  denote by  $\CI$ the
corresponding block; then one of the following holds:
\begin{enumerate}
\item   The  representation   $\varphi_\cs$  is   a  component   of  the
reflection  representation   $r$  of  $\WGL$.  In   this  case,  $\tilde
Q_\cs=q\widetilde\Id+\tphi_\cs$ and the block $\CI$ is regular.
\item  The  representation  $\varphi_\cs$   is  not  a  component  of
$r^{\wedge i}$  for any  $i$; then $\tQ_\cs=\tphi_\cs$. In
this case the block may or may not be regular.
\end{enumerate}
\end{proposition}
We  shall refer  to case  (i) by  saying that  $\cs$ is  {\it standard}.
Recall that a block  $\CI$ is regular if there exists  a local system in
$\CI$ with  support the regular class  and that in that  case this local
system  is unique  and  corresponds to  the  identity representation  of
$\WGL$ (\cf\  \cite[1.10]{DLM2}). \begin{proof} We prove  first that one
of the two properties for $\varphi_\cs$ and $\CI$ holds. This is done by
checking the  tables of the appendix.  First we reduce the  check to the
case  when  $\bG$  is  quasi-simple:  if $\bG$  is  not  quasi-simple  a
unipotent class is  a product of unipotent classes  of each quasi-simple
component and  a local system  on such a class  is the product  of local
systems  on the  components. In  particular  a subregular  class is  the
product of  the regular classes  of all the  components but one  and the
subregular class in the last component. A cuspidal datum is a product of
cuspidal data  for the quasi-simple  components. All this shows  that we
can reduce the verification to the quasi-simple case.

It  is  then  apparent  from  the  tables  that  when  $\varphi_\cs$  is
the   reflection  representation,   the  block   is  regular   and  that
otherwise  $\varphi$ has  dimension  strictly less  than the  reflection
representation,  so  appears in  no  exterior  power of  the  reflection
representation.

We  now prove  the formula  for  $\tQ_\cs$  in each  case. We  know
that $P_{\ci,\cj}$  is zero unless  $C_\ci\subsetneq\overline{C_\cj}$ or
$\ci=\cj$. So  $P_{\cs,\ci}=0$ unless  $C_\ci$ is  the regular  class or
$\ci=\cs$.

Consider first  the case when  $\CI$ is  regular: denote by  $\creg$ the
unique  pair in  $\CI$ with  regular support.  If we  take the  rows and
columns pertaining to  $\cs$ and $\creg$ to be the  last two, the matrix
equation
$\tilde P\inv\tilde \Lambda\inv(\lexp t{\tilde P}\inv)=\tilde\Omega$
which determines $\tilde P$ and $\tilde \Lambda$ has the form:
$$
\begin{pmatrix}\ldots&\ldots&\ldots\\0&1& Q\\0&0&1\end{pmatrix}
\begin{pmatrix}\ldots&\ldots&0\\0&\mu_\cs&0\\0&0&\mu_\creg\end{pmatrix}
\begin{pmatrix}\ldots&0&0\\\ldots&1&0\\\ldots&Q&1\end{pmatrix}=
\begin{pmatrix}\ldots&\ldots&\ldots\\\ldots&\tilde\omega_{\cs,\cs}&
\tilde\omega_{\cs,\creg}\\\ldots&\tilde\omega_{\creg,\cs}&
\tilde\omega_{\creg,\creg}\end{pmatrix}
$$
where $Q=(\tilde P\inv)_{\cs,\creg}$,
$\mu_\cs=(\tilde\Lambda\inv)_{\cs,\cs}$ and
$\mu_\creg=(\tilde\Lambda\inv)_{\creg,\creg}$.
We thus get:
$\mu_\cs+Q^2\mu_\creg=\tilde\omega_{\cs,\cs}$,
$Q\mu_\creg=\tilde\omega_{\cs,\creg}$ and
$\mu_\creg=\tilde\omega_{\creg,\creg}$.

In case (i) we apply \ref{tildeomegaik}. If $\bG_1,\ldots,\bG_k$ are the
quasi-simple components of $\bG$, we have
$r^{\wedge i}=\sum_{i_1+\ldots+i_k=i}r_1^{\wedge i_1}
\otimes\ldots\otimes r_k^{\wedge i_k}$,
where $r_i$ is the reflection representation of the $i$-th component of
$\WGL$.
So, using the remarks following \ref{preferred}
we have $\scal{\tphi_\cs}{\tilde r^{\wedge i}}\WGLF=
\begin{cases}1&\text{if $i=1$}\\0&\text{otherwise}\end{cases}$.
We then obtain $\tilde\omega_{\cs,\creg}=-
|Z^{0F}_\bG|q^{l-1}$ and $\tilde\omega_{\creg,\creg}=
|Z^{0F}_\bG|q^l$ where $l$ is as in \ref{tildeomegaik},
whence $Q=-q$, whence $\tilde P_{\cs,\creg}=q$.

In case (ii), the  above computation gives $\tilde\omega_{\cs,\creg}=0$,
so the only non-zero $\tilde P_{\cs,\ci}$ is $\tilde P_{\cs,\cs}=1$.

It remains  only to consider  case (ii)  for a non-regular  block, where
dimension  considerations imply  that the  only non-zero entry  $\tilde
P_{\cs,\ci}$ is $\tilde P_{\cs,\cs}$.

In either case, the value of $\tQ_\cs$ by is obtained by applying \ref{QiwF}.
\end{proof}
\begin{proposition}\label{restriction Q subregular}
Assume  that $\cs$  is an  $F$-stable  standard subregular  pair in  the
regular block $\CI_\bG$, and that $\bG$  is quasi-simple. Let $\bM$ be a
rational  Levi  subgroup  of  $\bG$, and  let  $C_1,\ldots,C_k$  be  the
$F$-stable subregular classes in $\bM$,  which are in bijection with the
set  of  $wF$-stable  irreducible  constituents $\bM_i$  of  $\Mo$.  Let
$\cs_i$ be  the pair corresponding  to the reflection  representation of
$W_{\bM_i}(\bL)$; then $\cs_i$ has support  $C_i$ and is a standard pair
in the regular block $\CIM$. Moreover we have
$$\Res^\WGLF_\WMLF\tQ_\cs=
\left((1-k)q\inv+\tphi_\cs(wF)-\sum_{i=1}^{i=k}\tphi_{\cs_i}(wF)\right)
\tQ_\cregM+\sum_{i=1}^{i=k}\tQ_{\cs_i}$$
where $\cregM$ is the pair with regular support in $\CIM$.
\end{proposition}
\begin{proof}
Let $V_\bG=Y(Z^0_\bL/Z^0_\bG)\otimes\bbR$, and
$V_\Mo=Y(Z^0_\bL/Z^0_\Mo)\otimes\bbR$. By \ref{preferred},
$\tphi_\cs$ is the extension of  the reflection representation of $\WGL$
which occurs in $V_\bG$, and by the same remarks we have $\Trace(vwF\mid
V_\Mo)=\sum_i\tphi_{\cs_i}(vwF)$ for  $v\in \WML$ (only  the $wF$-stable
components  occur when  we take  the trace  of an  element in  the coset
$\WMLF$). Thus if $V$ is the kernel of the natural map $V_\bG\to V_\Mo$,
we have $\Res^\WGLF_\WMLF\tphi_\cs=
\sum_{i=1}^{i=k}\tphi_{\cs_i}+\Trace(wF\mid V)\widetilde\Id$.
Evaluating both sides at $wF$ we get
$\Trace(wF\mid V)=\tphi_\cs(wF)-\sum_{i=1}^{i=k}\tphi_{\cs_i}(wF)$.

Now  by  \cite[1.10]{DLM2}  since  the block  $\CI_\bG$  is  regular  by
assumption, the  block $\CIM$ is also  regular. We know from  the remark
after  the statement  of \ref{Psr}  that the  pairs which  occur in  the
restriction of $\tQ_\cs$ have regular  or sub-regular support Since the
regular  class  corresponds to  $\widetilde\Id$  in  any regular  block,
$\cs_i$ must have support $C_i$, and thus $\cs_i$ is standard, so that by
\ref{Psr} we have $\tphi_{\cs_i}= \tQ_{\cs_i}-q\tQ_\cregM$.

The formula for the restriction of  $\tQ_\cs$ results from this and
the above formula for the restriction of $\tphi_\cs$.
\end{proof} From \ref{zeta2} and 
\ref{restriction Q subregular}
above, we deduce
\begin{proposition}\label{restriction Gamma subregular}
For any standard subregular pair $\cs$, we have
$$\sgn_\bG\sgn_\bM\sR^\bG_\bM\Gamma_\cs=\frac{a_\cs}
{a_{\cs_i}}\Gamma_{\cs_i}+\frac{a_\cs}{a_\cregM}\left((1-k)q+
\tphi_\cs(wF)-\sum_{i=1}^{i=k}\tphi_{\cs_i}(wF)\right)\Gamma_\cregM.$$
\end{proposition}

Similar computations can be made for non-standard pairs; however the end
result does not appear to have as clear a statement.
\section{The case of $\SL_n$}

We now  discuss the case  of $\bG=\SL_n$. According to  \cite[\S 5]{LS},
cuspidal data  are indexed by characters  of the centre $Z$  of $\SL_n$.
Assume that  $\chi$ is a  character of order  $d|n$ of $Z$;  then $\chi$
corresponds to an equivariant cuspidal local system on the regular class
of a Levi subgroup of type $A_{d-1}^{n/d}$. We will denote by $\CI_\chi$
the  corresponding  block  of  $\bG$. The  unipotent  classes  of  $\bG$
are  indexed  by  partitions  of  $n$.  Let  $C_\lambda$  be  the  class
indexed  by  the partition  $\lambda$  of  $n$.  There  is at  most  one
local  system  on  $C_\lambda$  in  $\CI_\chi$;  such  a  system  exists
when  all the  parts  of $\lambda$  are  divisible by  $d$  and we  will
denote it by  $\ci_\lambda^\chi$. When $\chi$ is  the trivial character,
$\ci_\lambda^\chi$ is the trivial local  system on $C_\lambda$, which is
also the  only irreducible  local system on  $C_\lambda$ in  $\GL_n$. We
will denote it simply by $\ci_\lambda$ in the latter case.

\begin{theorem}\label{Sln}  The   Laurent  polynomial  $\tilde
P_{\ci^\chi_\lambda,\ci^\chi_\mu}$ for  $\SL_n$ is equal to  the Laurent
polynomial  $\tilde  P_{\ci_{\lambda/d},\ci_{\mu/d}}$  for  $\GL_{n/d}$,
where $\lambda/d$ (resp. $\mu/d$) denotes  the partition whose parts are
$1/d$ times those of $\lambda$ (resp. $\mu$).
\end{theorem}
\begin{proof} The proof  consists of merely observing  that the equations
which determine $\tilde  P_{\ci^\chi_\lambda,\ci^\chi_\mu}$
and $\tilde P_{\ci_{\lambda/d},\ci_{\mu/d}}$ coincide.
In either case the equation may be  written:
$\tilde  P\inv\Lambda_1(\lexp t{\tilde P}\inv)=\Omega_1$ where
$\Lambda_1=|Z^{0F}_\bG|\inv\tilde\Lambda\inv$ and
$\Omega_1=|Z^{0F}_\bG|\inv\tilde\Omega$.
In the present case, $F$ acts  trivially   on   $\WGL$.   If,  for
$\varphi\in\Irr(\WGL)$,  we  denote  by  $\ci_\varphi$  the
corresponding  local system,  we have  according to  \ref{tildeomegaik}:
$$(\Omega_1)_{\ci_\varphi,\ci_{\varphi'}}=\sum_{i=0}^l q^{l-i}(-1)^i
\scal{\varphi\otimes\varphi'}{r^{\wedge  i}}\WGL.$$
We  have two  cases  to  consider: firstly  $\bG=\SL_n$,  $\bL$ of  type
$A_{d-1}^{n/d}$  and secondly  $\bG=\GL_{n/d}$, $\bL$  a maximal  torus.
In  either case  we  have $\WGL\simeq\gothS_{n/d}$  and $l=n/d-1$.  Thus
the  matrices $\Omega_1$  in the  two  cases may  be identified  through
the  bijection which  maps the  local system  $\ci^\chi_\lambda$ to  the
local system  $\ci_{\lambda/d}$ (since,  according to \cite[\S 5]{LS}
both  correspond under  the generalized  Springer correspondence  to the
character of  $\gothS_{n/d}$ indexed  by the partition  $\lambda/d$). To
verify that  the equations are the  same, it remains only  to check that
the rows and columns of the matrix $\tilde P$, 
both of which are
indexed by the irreducible characters  of $\gothS_{n/d}$, are ordered in
the same  way in either case.  This ordering is  induced by  the partial
order on unipotent  classes in either case, and  the coincidence follows
from the description of this partial order in terms of partitions: we
have  $C_\lambda\ge C_\mu$  if  and only  if  $\lambda\ge\mu$ where,  if
$\lambda=\{\lambda_1,\lambda_2,\ldots\}$  with $\lambda_1\ge\lambda_2\ge
\ldots$   (resp.  $\mu=\{\mu_1,\mu_2,\ldots\}$   with  $\mu_1\ge\mu_2\ge
\ldots$) this means that for all $i$ we have $\lambda_1+\ldots+\lambda_i
\ge\mu_1+\ldots+\mu_i$. This  condition is compatible with  dividing all
parts of $\lambda$ and $\mu$ by the same integer $d$, whence the result.
\end{proof}

The  significance of  the  previous  result is  that  in view of 
\ref{*RMG(Gamma)},  the
computation  of $\sR^\bG_\bM$  of the generalized Gelfand-Graev characters,
hence of the $\CX_\ci$, and through them of the $\CY_\ci$, and hence of the 
characteristic functions of the unipotent conjugacy classes for the group
$\SL_n$, is reduced to the same problem for various $\GL_{n'}$, which is
in principle known. According to the program in \cite{DLM1}, this is a
step towards determining the character table of $\SL_n(q)$. The other 
essential step in this program is the determination of $\sR^\bG_\bM$
of the irreducible characters, for which the work of C. Bonnaf\'e 
gives a solution.

\section{Appendix: local systems on the subregular unipotent class in
good characteristic for simply connected groups}
We  describe  now  the  generalized Springer  correspondence  for  local
systems  on  the  subregular  class for  simply  connected  quasi-simple
groups.  The  description  for  arbitrary  quasi-simple  groups  follows
easily.

This appendix contains information extracted from \cite{ICC}, \cite{LS} and
\cite{S}. The table below is as follows:
\begin{itemize}
\item The column ``$\bG$'' contains the type of $\bG$.
\item The column ``$C$'' describes the subregular class $C$, in Carter's
notation for exceptional  groups and by giving  the partition associated
to the Jordan form for classical groups.
\item  The column  ``Dynkin-Richardson'' contains  the Dynkin-Richardson
diagram of $C$.
\item The  column ``$A(u)$'' describes  the group $A(u)$ for  an element
$u\in C$.
\item The column ``$\ci$'' describes the local system $\zeta$ considered
on  $C$;  it is  described  by  giving  the  name of  the  corresponding
character  of $A(u)$;  this last  group  is when  possible described  as
a  Coxeter  group  so  the  naming  scheme  for  characters  of  Coxeter
groups  (see below)  applies. The  exceptions  are the  cyclic group  of
order 3  whose characters are  denoted $1,\zeta,\zeta^2$ and  the cyclic
group  group of  order 4  whose characters  are denoted  $1,i,-1,-i$. If
$\ci=(C,\zeta)$ let  $(\bL,\ci_0)$ be the corresponding  cuspidal datum,
where $\ci_0=(\zeta_0,C_0)$. In general there  is only one cuspidal pair
in $\bL$ (which is in most cases a local system on the regular class) so
neither $C_0$  nor $\zeta_0$  is mentioned; when  there is  an ambiguity
they are mentioned in the last column.
\item When  $\bL$ is not  a maximal torus $\bT$  or equal to  $\bG$, the
column ``$\bL$'' describes the Levi  by circling the nodes corresponding
to simple  roots of  $\bL$ on  the Dynkin diagram  of $\bG$.  The simple
roots of $\WGL$  in $X(Z^0_\bL/Z^0_\bG)\otimes\bbR$ therefore correspond
to the unmarked nodes of the same diagram.
\item  When  $\WGL$ is  neither  trivial  nor  equal  to $W_\bG$  it  is
described in the column ``$\WGL$'' by its Dynkin diagram, which has been
decorated  by  letters $a$,  $b$,  $\ldots$  which  appear also  on  the
un-circled nodes in the column  ``$\bL$'' to describe the correspondence
between simple reflections.
\item  The column  ``$\varphi_\ci$'' describes  the character  of $\WGL$
corresponding to $\ci$. The notation for characters of Coxeter groups is
as follows: $1$,  $\sgn$ and $r$ always represent the  trivial, sign and
reflection  representation  respectively.  Other linear  characters  are
represented  by  the  Dynkin  diagram  labelled by  the  values  of  the
character  on  the  simple  reflections.  The  notation  for  characters
of  $F_4$  is  that  from \cite{Carter}  (the  character  $\phi''_{2,4}$
factors through $W(F_4)/W(D_4)=W(A_2)$ and is trivial on the reflections
corresponding to  a short root;  the character $\phi'_{2,4}$  is deduced
from it  by the  diagram automorphism). The  characters of  $W(B_n)$ are
parametrized in the usual way by pairs of partitions.
\end{itemize}
\bigskip
\def\tbl#1{{\small
\offinterlineskip
\halign to \hsize
{\tabskip0pt\vrule##\vrule height10pt depth10pt width0pt\tabskip2ptplus12pt
 &\hfill$##$\hfill
 &\hfill$##$\hfill
 &\hfill$##$\hfill
 &\hfill$##$\hfill
 &\hfill$##$\hfill
 &\hfill$##$\hfill
 &\hfill$##$\hfill
 &\hfill$##$\hfill&\tabskip0pt\vrule##\cr
\noalign{\hrule}
&\bG&C&\text{Dynkin-Richardson}&A(u)&\ci&\bL&\WGL&\varphi_\ci&\cr
\noalign{\hrule}
#1
\noalign{\hrule}
}}}
\def\barlength{10 pt}
\def\node{{\kern -1pt\circ\kern -1pt}}
\def\bulnode{{\kern -1pt\bullet\kern -1pt}}
\def\dbulnode#1{{\kern-1pt\mathop\bullet\limits^{\hbox to 0pt{\hss
$\scriptstyle#1$\hss}}\kern-1pt}}
\def\rtbar{{\rlap{\vrule width\barlength height1.5pt depth-1pt}
           \rlap{\vrule width\barlength height2.5pt depth-2pt}
   \rlap{\hbox to\barlength{\hfill\raise 0.9pt\hbox{$\scriptstyle>$}\hfill}}
		 \vrule width\barlength height3.5pt depth-3pt}
		 }
\def\rdbar{{\rlap{\vrule width\barlength height2pt depth-1.5pt}
   \rlap{\hbox to\barlength{\hfill\raise 0.9pt\hbox{$\scriptstyle>$}\hfill}}
                 \vrule width\barlength height3pt depth-2.5pt}
		 }
\def\ldbar{{\rlap{\vrule width\barlength height2pt depth-1.5pt}
   \rlap{\hbox to\barlength{\hfill\raise 0.9pt\hbox{$\scriptstyle<$}\hfill}}
                 \vrule width\barlength height3pt depth-2.5pt}
		 }
\def\bar{{\vrule width\barlength height2.75pt depth-2.25pt}}
\def\vertbar#1#2{\rlap{\kern 0.7pt\vrule width0.5pt height1pt depth9pt}
		 \rlap{\lower12.7pt\rlap{$#2$}}#1}
\def\sm{{\hbox{-}}}
\def\below#1#2{\vtop{\hbox{\strut#1}\hbox{\strut#2}}}
\tbl{
&G_2&G_2(a_1)&\dbulnode2\rtbar\dbulnode0&W(A_2)&
   1&\bT&W(G_2)&r&\cr
&&&&&r&\bT&W(G_2)&\dbulnode{-1}\rtbar\dbulnode1&\cr
&&&&&\sgn&\bG&1&1&\cr
&F_4&F_4(a_1)&\dbulnode2\bar\dbulnode2\rdbar\dbulnode0\bar\dbulnode2&W(A_1)&
    1&\bT&W(F_4)&r&\cr
&&&&&\sgn&\bT&W(F_4)&\phi'_{2,4}&\cr
&E_6&E_6(a_1)&\dbulnode2\bar\dbulnode2\bar
\vertbar{\dbulnode0}{\bulnode\scriptstyle2}\bar\dbulnode2\bar\dbulnode2
&\bbZ/3\bbZ&
   1&\bT&W(E_6)&r&\cr
&&&&&\zeta&\node\bar\node\bar\vertbar{\dbulnode a}{\bulnode\scriptstyle b}
\bar\node\bar\node&\dbulnode b\rtbar
\dbulnode a&\dbulnode{-1}\rtbar\dbulnode 1&\cr
&&&&&\zeta^2&\multispan 3{\vtop{\parindent=0pt
\hsize=19em
same description. The cuspidal local system
is the other one on the regular class of
$\bL\simeq \SL_3\times_{Z(\SL_3)}\SL_3$}}&\cr
}\vfill\eject\tbl{
&E_7&E_7(a_1)&\dbulnode2\bar\dbulnode2\bar\dbulnode2\bar\vertbar{\dbulnode0}
{\bulnode\scriptstyle2}\bar
\dbulnode2\bar\dbulnode2&W(A_1)&
   1&\bT&W(E_7)&r&\cr
&&&&&\sgn&\node\bar\dbulnode a\bar\node\bar\vertbar{\dbulnode b}\node
\bar\dbulnode c\bar
\dbulnode d&\dbulnode a\bar\dbulnode b\ldbar\dbulnode c\bar\dbulnode d&
\phi''_{2,4}&\cr
&E_8&E_8(a_1)&\dbulnode2\bar\dbulnode2\bar\dbulnode2\bar\dbulnode2\bar
\vertbar{\dbulnode0}{\bulnode\scriptstyle2}\bar\dbulnode2\bar\dbulnode2
&1&1&\bT&W(E_8)&r&\cr
&\below{$A_n$}{$n$ even}&(1,n\sm 1)&
\displaystyle\dbulnode2\cdots\dbulnode2\bar\dbulnode0\bar\dbulnode2
\cdots\dbulnode2
&1&1&\bT&W(A_n)&r&\cr
&\below{\strut$A_n$}{$n$ odd}&(1,n\sm 1)&
\dbulnode2\cdots\dbulnode2\bar\dbulnode1\bar\dbulnode1\bar\dbulnode2
\cdots\dbulnode2
&1&1&\bT&W(A_n)&r&\cr
&B_n&(1,1,2n\sm 1)&\dbulnode2\cdots\dbulnode2\rdbar\dbulnode0&
W(A_1)&1&\bT&W(B_n)&r&\cr
&&&&&\sgn&\bT&W(B_n)&(1.n\sm 1,\emptyset)&\cr
&C_2&(2,2)&\dbulnode0\ldbar\dbulnode2&W(A_1)&1&\bT&W(C_2)&r&\cr
&&&&&\sgn&\bT&W(C_2)&(\emptyset,2)&\cr
&\below{$C_n$}{$n>2$}&(2,2n\sm 2)
&\dbulnode2\cdots\dbulnode2\bar\dbulnode0\ldbar\dbulnode2&
W(A_1)^2&(1,1)&\bT&W(C_n)&r&\cr
&&&&&(\sgn,\sgn)&\bT&W(C_n)&(\emptyset,n)&\cr
&&&&&(\sgn,1)&
\bulnode\cdots\bulnode\bar\bulnode\bar\bulnode\ldbar\node
&W(C_{n-1})&(\emptyset,1.n\sm 2)&\cr
&&&&&(1,\sgn)&
\bulnode\cdots\bulnode\bar\node\bar\node\ldbar\node
&W(C_{n-3})&1&\cr
&\below{$D_n$}{$n$ odd}&(3,2n\sm 3)&
\dbulnode 2\cdots\dbulnode2\bar\vertbar{\dbulnode 0}
{\bulnode\scriptstyle 2}\bar\dbulnode2
&\bbZ/4\bbZ&1&\bT&W(D_n)&r&\cr
&&&&&-1&
\bulnode\cdots\bulnode\bar\node\bar\vertbar\node\node\bar\node
&W(B_{n-2})&(1.n\sm 3,\emptyset)&\cr
&&&&&i&
\below{$\node\bar\bulnode\cdots\node\bar\bulnode\bar\node\bar\node\bar\vertbar
\node\node\bar\node$}{(type $D_5\times A_1^{(n-5)/2}$)}
&W(B_{\frac{n-5}2})&1&\cr
&&&&&-i&\multispan 3{\vtop{\parindent=0pt\hsize=19em
same description. The cuspidal local system
is also parametrized by $-i$ on the $D_5$ component of $\bL$
}}&\cr
&\below{$D_n$}{$n$ even}&(3,2n\sm 3)&
\dbulnode 2\cdots\dbulnode2\bar\vertbar{\dbulnode 0}
{\bulnode\scriptstyle 2}\bar\dbulnode2
&W(A_1)^2&(1,1)&\bT&W(D_n)&r&\cr
&&&&&(-1,1)&
\bulnode\cdots\bulnode\bar\node\bar\vertbar\node\node\bar\node
&W(B_{n-2})&(1.n\sm 3,\emptyset)&\cr
&&&&&(1,-1)&
\below{$\node\bar\bulnode\cdots\node\bar\bulnode\bar\node\bar\vertbar\bulnode
\bulnode\bar\node$}{(type $A_1^{n/2}$)}
&W(B_{n/2})&(\emptyset,n/2)&\cr
&&&&&(-1,-1)&
\below{$\node\bar\bulnode\cdots\node\bar\bulnode\bar\node\bar\vertbar\bulnode
\node\bar\bulnode$}{(type $A_1^{n/2}$)}
&W(B_{n/2})&(\emptyset,n/2)&\cr
}

\medskip
A more precise description of the local systems when
$\bG=\text{Spin}_{2n}$ (the simply connected
semi-simple group of type $D_n$) is as follows:
$A_{\text{SO}_{2n}}(u)$ is isomorphic to $W(A_1)$;
when $n$ is odd it is the unique subgroup of order 2 of $A_\bG(u)$, while
when $n$ is even it is the first $W(A_1)$ in $A_\bG(u)$.

\end{document}